\magnification=\magstep1
\input mssymb

\overfullrule=0pt
\def\sqr#1#2{{\vcenter{\hrule height.#2pt
      \hbox{\vrule width.#2pt height#1pt \kern#1pt
         \vrule width.#2pt}
       \hrule height.#2pt}}}
\def\op{\mathchoice\sqr34\sqr34\sqr{2.1}3\sqr{1.5}3}

\def\Jbeop{$\op $}
\def\eop{\ifmmode\op\else\Jbeop\fi }
\def\NUOLI#1#2{\setbox0=\hbox to #2truecm{\rightarrowfill}
\mathop{\box0}\limits_{#1}} 
\newdimen\refskip  
\refskip=35pt                                         
\def\ref #1 #2\par{\par\noindent\rlap{#1}\hskip\refskip
\hangindent\refskip #2}
\def\kohta #1 #2\par{\par\indent\rlap{\rm#1}\hskip\refskip
\hangindent40pt #2\par}
\null
\vskip 3truecm
\centerline{\bf Structure of local Banach spaces of locally convex spaces}
\footnote{}{Mathematics subject classification: primary 46A05, secondary
46B20, 47A68 }
\bigskip\bigskip\bigskip
\centerline{Jari Taskinen, University of Helsinki}
\bigskip\bigskip\bigskip\bigskip
{\bf 1. Introduction and notation.} In this paper the word "local" is
used in at least three different meanings.
Our aim is to  study local Banach spaces
of Fr\'echet or other locally convex spaces, and it turns out that it
is convenient to use the local theory of Banach spaces for this purpose.
Recall that given a locally convex
space  $E$  and a continuous seminorm  $p$  on  $E$ the completion
of the normed space  $(E/{\rm ker}(p),p)$  is a local Banach space of
$E$, and it is denoted by  $E_p$.  If  $(p_{\alpha})_{\alpha \in A}$  is
a system of seminorms defining the topology of  $E$, then we call
$(E_{p_{\alpha}})_{\alpha \in A}$  a system of local Banach spaces.

Valdivia proved in [V] that every infinite dimensional
nuclear locally convex space  $E$  has
for all separable infinite dimensional Banach spaces
$X$  a system of local Banach spaces
isomorphic to  $X$.  In this paper we study Schwartz spaces and,
consequently, compact linear operators in Banach spaces. We consider
''dense'' factorizations of operators in Banach spaces,
i.e., factorizations $T = T^{(2)} T^{(1)}$,
$T \in L(X,Y),$ through $Z$ such that $T^{(1)} (X) \subset Z $ is dense.
Assuming density is enough to guarantee a direct application to the
study of local Banach spaces of locally convex spaces.

In Section 2. we introduce a way to
measure the distance of nonisomorphic Banach spaces, the so called local
distance function.  We prove  general theorems (Theorems 2.8 and 2.14)
which show
that a compact operator  $T$  on a separable reflexive Banach space  $X$
factors through a Banach space  $Y$  in a strong sense (the image of
$X$ is dense in $Y$; both factors of the given operator are compact),
if the distance
of  $X^*$  and  $Y^*$  is small enough and some technical assumptions are
satisfied. For example we show that an arbitrary compact operator on a
separable ${\cal L}_p$--space $X$, $1 < p < \infty$, factors through any
separable ${\cal L}_p$--space $Y$ in the above strong sense (Theorem
2.10).
Consequently, if $X$ and $Y$ are such spaces, we see that a Schwartz
space having a system of local Banach space isomorphic to $X$ also has
a system isomorphic to $Y$ (Corollary 2.16).

The concept of local distance function is also analyzed from the Banach
space theoretical point of view.  We show (Proposition 2.4) that there exist
separable, reflexive Banach spaces which are in a sense very close to
$\ell_2$  but which do not have a Schauder basis.  This result is
a consequence of a construction of Szarek, [S].

Section 3 contains some more remarks on dense factorizations. The
considerations in this section are technically easy, compared to Section
2, and the results are not so deep from the Banach--space theoretical
point of view. But the consequences to the locally convex space theory
are so strong that it is worthwhile to present the results in detail.
In Theorem 3.2 we show that if $X $ is a separable Banach space with a
complemented unconditional basic sequence, if $Y$ and $Z$ are Banach
spaces, $Z$ separable, and if $T \in L(X,Y)$ is compact, then $T$
factors through $X \times Z$ such that the image of $X$ is dense in
$X \times Z$. Combining this observation with the results of  Section 2
we see that a compact $T \in L ( {\ell }_2)$, $1 < q < \infty$, factors
densely  through any separable ${\cal L}_p$--space for $1<p \leq
\infty$.
Another
consequence of Theorem 3.2 is that a compact $T \in L(X)$ factors
through $Y \times Z$ for an arbitrary separable Banach space $Z$, if
the (uniform, see later) local distance of $X^*$ and $Y^*$ is small
enough and if some technical assumptions are satisfied.

In the second part of the paper, Section 4, we study some duality problems for local
Banach spaces. The relations of the local Banach spaces of a locally
convex space $E$ and of the strong dual $E'_b$ are not yet well
understood in general. For example,
it is an open problem whether, given a Fr\'echet or a
$(DF)$--space  $E$  with a system of local Banach spaces isomorphic to
some Banach space  $X$,  $E$ also  has a fundamental system of
Banach discs $(B_{\alpha})_{\alpha \in A}$ such that the corresponding
Banach spaces $E_{B_{\alpha}}$ are isomorphic to $X.$
We are not able to solve this
problem here, but in Proposition 4.1 we generalize the known partial
positive results, and
in Proposition 4.6 we construct an example  which can be
considered as a partial negative solution.  This construction uses
an estimate (Lemma 4.2) of the absolute projection constant of a Banach
space.  The method used in the proof of this lemma may be new.

We use mainly the terminology of [K1], [LT1] and [TJ].  Let us however
mention some notations and definitions.  ${\Bbb N}$  stands for the
set  $\{1,2,3,\ldots\}$.  The closure of a set $A$ is denoted by
$\overline{A}$. A subset $A$ of a vector space is absolutely convex,
if $\sum_{i=1}^n \lambda_i x_i \in A$ for all sequences
$(x_i)_{i=1}^n \subset A$
and for all scalar sequences $(\lambda_i)$  satisfying
$\sum_{i=1}^n | \lambda_i | \leq 1.$ A closed, bounded and absolutely
convex subset of a Fr\'echet space is called a Banach disc.
The vector spaces are over the real or
complex scalar field unless otherwise stated.  By a subspace we mean
a linear subspace, by an operator a continuous linear operator, and by
an isomorphism a linear homeomorphism.  Two Banach spaces  $X$  and  $Y$
are $C$--isomorphic, $C \geq 1$, if the Banach--Mazur distance
$d(X,Y) \leq C$, where  $d(X,Y) := \inf\ \Vert \psi \Vert \Vert
\psi^{-1}\Vert,\ \psi : X \to Y$  isomorphism.  A subspace  $X$  of
a Banach space  $Y$  is $C$--complemented, if there exists a projection
$P$  from  $Y$  onto  $X$  with  $\Vert P\Vert \leq C$; the
projection constant  $ \lambda  (X,Y) $  is the number  $\inf \ \{\Vert
P\Vert \mid P$  is a projection from  $X$  onto  $Y\}$.  The absolute
projection constant $\lambda (X)$ of a Banach space  $X$  is the supremum of
$\lambda (Z,Y)$  over all Banach spaces  $Y$  containing a subspace
$Z$  isometric to  $X$.

For the definition of the bounded approximation property, (Schauder)
basis and basis constant we refer to [LT1].  The definition and
properties of ${\cal L}_p$--spaces can be found in [LT],
Chapter 5, or in [TJ], \S 2.
Type, cotype and the corresponding constants are defined in [TJ],
\S 4.  Some elementary facts on and definitions of tensor products,
especially the projective tensor product, are mentioned in [TJ], \S 5
(or in [T]; for more details, see [J]),
and the real interpolation method is also presented in the
same book, \S 3. 

Most of these definitions occur also in [J].

If  $X$  and  $Y$  are Banach spaces we denote by  $L(X,Y)$  the space
of bounded linear operators  $X \to Y$.  If  $T \in L(X,Y)$ is
compact, we define
the $n$:th approximation number $a_n(T),$ $ \ n \in {\Bbb N}$,  by
$$a_n(T) := \inf\ \{\Vert T - T_n\Vert \mid T_n \in L(X,Y),\ {\rm rank}\
T_n < n\}.\eqno(1.1)$$
\indent Recall that a locally convex space  $E$  is a Schwartz space,
if for every continuous seminorm  $p$  there exists a continuous
seminorm  $q,\ q \geq p$, such that the canonical mapping  $E_q \to E_p$
induced by the identity operator on  $E$  is compact.  For other
definitions concerning locally convex spaces we refer to [K1] and [J].

{\bf Acknowledgements.} I am very grateful to Jos\'e Bonet for reading
a preliminary version of the manuscript and for valuable remarks.
I would like to thank Hans--Olav Tylli for many
discussions on the subjects of this paper. I would also like to thank
Mikael Lindstr\"om for remarks and comments.
\bigskip
{\bf 2. Local distance function of Banach spaces with applications to
compact operators and Schwartz spaces.}  The Banach--Mazur distance
provides a natural way to measure differences of two isomorphic Banach
spaces, but for nonisomorphic spaces the Banach--Mazur distance is not
a finite number.  In this section we introduce a way to measure the
distance of nonisomorphic spaces, the so called local distance
function.  The local distance function is not a metric,
since the ''local distance'' of two Banach spaces is a function, not a
number. Some basic properties of local distance function are
given in Remarks 1.--8.
\bigskip
{\bf 2.1. Definition.} Let  $X$  and  $Y$  be Banach spaces and  $C \geq 1$.
The $C$--local distance function of  $X$  and  $Y$  is the function
$(C){\rm -}ld(X,Y): {\Bbb N} \to {\Bbb R}^+ \cup \{ \infty \} $,
the value of which at $n \in {\Bbb N}$ is the
infimum of ${\infty}$ and all numbers $f(n) \in {\Bbb R}^+  $
satisfying the following: If $M \subset X$  and  $M_Y \subset Y$  are
at most
$n$--dimensional subspaces, $n \in {\Bbb N}$, then there exist
$C$--complemented finite dimensional subspaces  $N \subset X$  and
$N_Y \subset Y$, such that  $M \subset N,\ M_Y \subset N_Y,\
{\rm dim}\ (N) = {\rm dim}\ (N_Y)$  and
$$d(N,N_Y) \leq f(n).$$
\indent Let  $X,Y$  and  $C$  be as above and let  $K: {\Bbb N} \to
{\Bbb N}$  be a non--decreasing function.  The $K$--uniform $C$--local
distance function of  $X$  and  $Y$  is the function  $(K,C){\rm -}ld(X,Y):
{\Bbb N} \to  {\Bbb R}^+  \cup \{ \infty \} $ the value of which
at the point $n \in {\Bbb N}$ is the infimum of ${\infty}$ and all
numbers $f(n) \in {\Bbb R}^+ $
satisfying the
following condition: If  $M \subset X$  and  $M_Y \subset Y$  are
subspaces with dimension not greater than
$n$, then there exist $C$--complemented subspaces
$N \subset X$  and  $N_Y \subset Y$  containing  $M$  and  $M_Y$,
respectively, such that  ${\rm dim}\ (N) = {\rm dim}\ (N_Y) \leq K(n)$
and  $d(N,N_Y) \leq f(n)$.
\bigskip
We call the $C$--local distance function just the local distance
function, if $C$ is clear from context or does not need to become
specified. In the same way we speak about the uniform local distance
function.

We say that the local distance (respectively, uniform local distance)
of  $X$  and  $Y$  is bounded, if  $(C){\rm -}ld(X,Y)$  is a bounded function
for some   $C$  (resp. $(K,C){\rm -}ld(X,Y)$  is a bounded function for some
$C$  and  $K$).

We say that the local distance of  $X$  and  $Y$  is finite, if, for
some  $C \geq 1,\  (C){\rm -}ld(X,Y)(n)$  is finite for all  $n$.

In the same way we define what means that the uniform local distance
of  $X$  and  $Y$  is finite.
\bigskip
{\bf 2.2. Remarks.} 1. It may happen that the local distance of  $X$
from itself is not finite.  This is the case for the Banach space  $X$
constructed by Pisier, [Ps1]: it has the property that all
$n$--dimensional subspaces are not better than $C\sqrt n$--complemented
in  $X$, where  $C$ is a constant depending on the space $X$ only.
In fact, the property
that the local distance of  $X$  from itself is finite is equivalent to
$X$  being a $\pi$--space in the sense of [JRZ].  This property is
stronger than the bounded approximation property.  We refer to [JRZ],
Proposition 1.1 and 1.2.

2. It would also be possible to define another
local distance function roughly speaking as the infimum of all functions
$f(n)$  satisfying the following condition ($X,$ $Y$ and $C$ as in
Definition 2.1): Given finite dimensional subspaces $M \subset X,$
$M_Y \subset Y$ there exist $C$--complemented subspaces $N$ and $N_Y$
containing $M$ and $M_Y,$ respectively, such that dim$(N)=$ dim$(N_Y)$
and such that $d(N,N_Y) \leq f($dim$(N)).$  We
would get a distance function which would be easier to estimate for
example in the case of $L_p$--spaces (cf.\ Remark 6.).  However, this
function seems not to be so useful in considerations like Theorem 2.8.

Omitting the space  $M_Y$  in Definition 2.1 would lead to another
concept of distance function.  This distance function would not be
symmetric with respect to  $X$  and  $Y$.

3. If $X$ and $Y$ are ${\cal L}_p$--spaces, then  the uniform
local distance of  $X$  and  $Y$
is bounded. This statement follows immediately
from definitions and [PR], Corollary 2.1. We do not know, if
the converse statement (''if $X$ is a ${\cal L}_p$-space
and the local distance of $X$ and $Y$ is bounded, then $Y$ is
a ${\cal L}_p$-space'') also holds.

4. Let  $X$  (resp.  $Y$)  be a Banach space of type  $p$  (resp.
$p'$).  If the local distance of $X$ and $Y$ is bounded,
then  $p = p'$.  Suppose by antithesis that  $1 \leq p < p' \leq 2$.
By the Maurey--Pisier theorem, [TJ], Theorem 7.6, $X$  contains a
subspace  $M$  which is 2--isomorphic to   $\ell_p^k,\ \infty > k >
(4C'D)^{pp' / (p-p')}$, where  $D = \sup\limits_n (C){\rm -}ld(X,Y)(n)$ and
$C' = T_{p'}(Y) $, the type $p'$--constant of  $Y$.

Let  $N \subset Y$  be any subspace which is $D$--isomorphic to a subspace
of  $X$  containing   $M$.  Then we have for the type $p'$--constant of
$N$  (see [TJ], (4.5), p.\ 15),
$$T_{p'}(N) \geq D^{-1}T_{p'}(M) \geq D^{-1}k^{1/p'-1/p}/2 \geq 2C',$$
which contradicts  $T_{p'}(Y) = C'$.

In the same way one shows that if the local distance of  $X$
and  $Y$  is bounded, then the spaces  $X$  and  $Y$  have the same
cotype.

5. It is a direct consequence of definitions that if  $ (C){\rm -}ld(X,Y)$
is a bounded function for some  $C$, then  $X$  (resp.  $Y$)  is crudely
finitely representable (and even strongly representable)
in  $Y$  (resp.  $X$). For the definition of these
concepts, see [MS], 11.6 and [BDG].

6. We show that if  $X$  (resp.  $Y$) is an $L_p$--space (resp.
$L_{p'}$--space), $1 \leq p \leq p' \leq \infty$, then, for all
$\varepsilon > 0$, for all  $n \in {\Bbb N}$,
$$(K,C){\rm -}ld(X,Y)(n) \leq
(1+\varepsilon ) (2(n+1)^2 / C(\varepsilon))^{n/p-n/p'}, \eqno(2.1)$$
where $C = 1+\varepsilon$ and $K(n) =
(2(n+1)^2 / C(\varepsilon))^n $ and $C(\varepsilon)$ is a constant
depending on $\varepsilon$.
Pe{\l}czy\'nski, Rosenthal and Kwapien have shown ([PR],
Theorem 2.1) that, given an $n$--dimensional subspace  $M$  of an
$L_q$--space, there exists a $1 + \varepsilon$--complemented subspace
$N \supset M$  such that  ${\rm dim}\ (N) \leq (2(n+1)^2 / C(\varepsilon))^n  $
and  \hbox{$d(N,\ell_q^{{\rm dim}\ (N)}) \leq
1+\varepsilon /3 $.}  So,
given $n$--dimensional subspaces  $M \subset X$  and  $M_Y \subset Y$,
we choose $k$--dimensional $1 + \varepsilon$--complemented subspaces
$N \supset M,\ N_Y \supset M_Y$  such that
$k \leq (2(n+1)^2 / C(\varepsilon))^n,$
$d(N,\ell_p^k) \leq 1+\varepsilon /3 ,\ d(N_Y,\ell^k_{p'}) \leq 1+\varepsilon
/3 $.
Then we have
$$d(N,N_Y) \leq (1 + \varepsilon) k^{1/p-1/p'} \leq ( 1+\varepsilon )
(2(n+1)^2 / C(\varepsilon))^{n/p - n/p'} .$$

We do not know how sharp this estimate is.

7. We do not know what is the relation of the local distance function
of Banach spaces $X$ and $Y$ and on the other hand of the local
distance function of the duals $X^*$ and $Y^*$. This question
should be compared with [LT], II.5.7 and 8, and [M]. The local distance function
of $X$ and $Y$ gives direct information only on some quotients of
$X^*$ and $Y^*$.

8. Let $X$ be a weak Hilbert space in the sense of [Ps2]. If $M \subset
X$ is an $n$--dimensional subspace, then it is known that for a
constant $c$ depending on the space $X$ only we have $d(M,{\ell}_2^n ) <
c  \log (n+1);$ see
[Ps2], Corollary 2.5. According to some yet unpublished information
Maurey has improved the result of Johnson and Pisier, [JP], showing
that weak Hilbert spaces have the "linear uniform projection
property". This means that given $\varepsilon > 0 $ and
a finite--dimensional subspace $M$ of $X$
we can find a $1 + \varepsilon$--complemented subspace $N$ containing
$M $ such that dim$(N) \leq  c_1$dim$(M)$,
where $c_1 > 0$ depends on the space $X$ and $\varepsilon$ only.

Summing up the above statements we get for every $\varepsilon > 0$
$$ (K,C){\rm -}ld (X,{\ell}_2) (n) \leq c_0 \log(n+1),  $$
where $C = 1+ \varepsilon , $ $ K : {\Bbb N} \to {\Bbb N}$
is an asymptotically linear function depending
on $\varepsilon$ and $X,$ and the constant $c_0$ depends on the space $X$
and $\varepsilon$ only.
\bigskip
Proposition 2.4 below shows that the properties of Banach spaces
with a ''small'' (but not bounded) local distance may be quite
different.  For this result we need to use a construction of Szarek.
We begin with
\bigskip
{\bf 2.3. Lemma.} {\sl For all  $n \in {\Bbb N}$  and all  $q,\ 2 < q <
\infty$, there exists an $n$--dimensional subspace  $Y^n_q $ of
the ${\Bbb R}$--Banach space $ L_q$
such that
$$bc(Y^n_q \bigoplus_{\ell_2} F) \geq cn^{(1/2)(1/2-1/q)}D^{-1/2},\eqno(2.2)$$
for an absolute constant  $c$  and for all normed spaces  $F$  satisfying
the following: if  $M \subset F$  is an $n$--dimensional subspace, then
$d(M , \ell^n_2) \leq D$.}
\bigskip
Here $bc(X)$ denotes the basis constant of the Banach--space $X$.
Recall that $X$ has a basis if and only if $bc(X)$ is finite.

This lemma is an important result of Szarek, see [S], Proposition 3.1.
Now we can prove
\bigskip
{\bf 2.4. Proposition.}  {\sl Let  $f: {\Bbb N} \to {\Bbb R}^+$  be a
non--decreasing, unbounded function, $f(1) \geq 3$.  For all  $\varepsilon,\ 0 < \varepsilon < 1$,
there exists a separable reflexive
${\Bbb R}$--Banach space  $X$  having no Schauder basis such that,
for all  $n \in {\Bbb N}$,
$$(K,C){\rm -}ld(\ell_2,X)(n) \leq f(n),
\eqno(2.3)$$
where  $C = 1 + \varepsilon$
and  $K$  is some function}  ${\Bbb N} \to {\Bbb N}$.
\bigskip
The  function $K$ is specified in (2.13) below.

Since $X$ is reflexive, also the dual $X^*$ does not have a basis.
\bigskip
{\bf Proof.} We choose the sequences  $(n_k)^{\infty}_{k=1},\ n_k \in
{\Bbb N}$, and  $(q_k)^{\infty}_{k=1},\ 2 < q_k < \infty$, as follows:
Let  $q_0 = 4,\ n_0 = 1$, and assume that  $k \in {\Bbb N}$  and that
$q_t,n_t$  are chosen for  $t < k$.  We define  $q'_k,2 < q'_k < q_{k-1}$,
such that
$$n_{k-1}^{1/2-1/q'_k} \leq 1 + \varepsilon/3 \eqno(2.4)$$
and then  $n_k > 2 n_{k-1}$  such that
$$f(n_k) \geq kf(n_{k-1}), \ \ \ n_k^{1/2-1/q'_k} \geq f(n_{k-1})/2
\eqno(2.5)$$
and, finally, $q_k,\ 2 < q_k \leq q'_k$  such that
$$n_k^{1/2-1/q_k} = f(n_{k-1})/2.\eqno(2.6)$$

A useful refinement of this choice is described in Remark 2.5 below.

We define
$$X = (\bigoplus^{\infty}_{k=1} Y_k)_{\ell_2}, \eqno(2.7)$$
where  $Y_k := Y^{n_k}_{q_k}$  and the spaces  $Y^{n_k}_{q_k}$  are as
in Lemma 2.3.  We first show that (2.3) holds.  So let  $n \in {\Bbb N}$
and let  $M \subset X,\ {\rm dim}\ (M) \leq n$.  Let  $k \in {\Bbb N}$  be
such that  $n_{k-1} \leq n \leq n_k$.  We define the finite dimensional
subspace  $N$  by
$$N := \bigoplus_{t\leq k}\ Y_t \bigoplus Q(M), \eqno(2.8)$$
where  $Q$  is the natural projection from  $X$
onto  $\overline{ \bigoplus\limits_{t > k}\ Y_t }$.
Clearly, $M \subset N.$ We show that  $N$  is $1 +
\varepsilon$--complemented in  $X$. Let $Q_t$  be the natural projection
from  $X$  onto  $Y_t$. Each space
$Q_t(M)$  is at most $n_k$--dimensional.  Since every  $Y_t$  is, by
definition, a subspace of  $L_{q_t}(0,1)$, we can use a result of Lewis,
[L], Corollary 4, to find for each  $t > k$  a projection  $P_t$  from
$Y_t$  onto  $Q_t(M)$  such that
$$\Vert P_t\Vert \leq n_k^{1/2-1/q_t} \leq 1 + \varepsilon/3.  \eqno(2.9) $$
Hence, there exists a projection  $\tilde Q$  from
$\overline{\bigoplus_{t>k}\
Y_t } \subset X $  onto
$\overline{\bigoplus_{t > k }\ Q_t(M) }$  such that
$\Vert \tilde Q\Vert \leq 1 + \varepsilon/3$.  On the other hand, since  ${\rm dim}\
(Q_t(M)) \leq n_k$  for all  $t$, we get, for  $t > k$,
$$d(Q_t(M),\ell^{{\rm dim } \ ( Q_t(M) )}_2 )
\leq n_k^{1/2-1/q_t} \leq 1 + \varepsilon/ 3 ,  \eqno(2.10)$$
see [L], Corollary 5.  This means that
$$d( \overline{\bigoplus_{t > k }\ Q_t(M) } ,\ell^m_2) \leq 1 +
\varepsilon/ 3  \eqno(2.11)$$
where  $m = {\rm dim}\
(\overline{ \bigoplus\limits_{t > k}\ Q_t(M) })
\in {\Bbb N} \cup \{ \infty \} $.  Using the
orthonormal projection in Hilbert space we thus find a projection
$\hat Q$  from
$\overline{ \bigoplus\limits_{t > k}\ Q_t(M) }$  onto  $Q(M)$  with
$\Vert \hat Q\Vert \leq 1 + \varepsilon/3 $.
Now  $\hat Q \tilde Q$  is a projection from
$ \overline{ \bigoplus\limits_{t>k}\ Y_t } $  onto  $Q(M)$  satisfying
$\Vert \hat Q \tilde Q \Vert \leq 1 + \varepsilon $, and this means that there exists a
projection  $R$  from  $X$  onto  $N$
with  $\Vert R\Vert \leq 1 + \varepsilon$.

We have
$$\eqalignno{d((\bigoplus_{t\leq k}\ Y_t)_{\ell_2},\ell_2^{(n_1+\ldots
+n_k)}) &\leq \sup\limits_{t\leq k}\ d(Y_t,\ell_2^{n_t}) \leq
\sup\limits_{t\leq k}\ n_t^{1/2-1/q_t}\cr
\noalign{\vskip 6pt}
&\leq f(n_{k-1})/2 \leq f(n) ,&(2.12)\cr}$$
see (2.6).  So, (2.11) and (2.12) imply  $d(N,\ell_2^{{\rm dim}\ (N)})
\leq f(n)$.   This proves the statement
$$(C){\rm -}ld(X,\ell_2)(n) \leq f(n)$$
for all  $n$.

From (2.8) we see that  $ n_k \leq
{\rm dim}\ (N) \leq 3n_k$;
so, this upper estimate for dim($N$)  depends only on  ${\rm dim}\ (M)$
(see the choice of $k$) and we get
$$n_k \leq K(n) \leq  3 n_k   \eqno(2.13)$$
for $n_{k-1} \leq n \leq n_k$.
This proves the
uniform version of the estimate (2.3).

The proof of [S], Proposition 4.1, with small modifications, shows that
$X$  does not have even a ''local basis structure''.  In fact, to prove
that  $X$  does not have a basis, it is sufficient to show that for all
$t \in {\Bbb N}$ large enough, the $n_{t}$--dimensional subspaces $Z$ of
$(\bigoplus\limits_{k\neq t}\ Y_k)_{\ell_2}$  satisfy
$d(Z, \ell_2^{{\rm dim}(Z) } ) \leq f(n_{t-2})/2$ (see Lemma 2.3);
then we have, by (2.6) and (2.5)
$$\eqalign{bc(X) &= bc(Y_{t}\oplus \ (\bigoplus_{k\neq t}\ Y_k)_{\ell_2})\cr
\noalign{\vskip 6pt}
&\geq \sqrt {2} c\ n_{t}^{(1/2)(1/2-1/q_{t} )}  f(n_{t-2})^{-1/2}\cr
\noalign{\vskip 6pt}
&= cf(n_{t-1})^{1/2}f(n_{t-2})^{-1/2} \geq ct^{1/2}
\NUOLI{t\to \infty}{1} \infty.\cr}$$
Given  $M \subset (\bigoplus\limits_{k\neq t}\ Y_k)_{\ell_2},\ {\rm dim}\
(M) = n_{t}$, we have  $M \subset (\bigoplus\limits_{k\neq t}
Q_k(M))_{\ell_2}$.  We then get, analogously to (2.10) and (2.12),
$$d(Q_k(M),\ell_2^{{\rm dim}\ (Q_k(M))}) \leq 1 + \varepsilon$$
for  $k > t$, and
$$d(Q_k(M),\ell_2^{{\rm dim}\ (Q_k(M))}) \leq n_{t-1}^{1/2-1/q_{t-1}}
= f(n_{t-2})/2$$
for  $k < t$.  This implies
$$d(M,\ell_2^{{\rm dim}\ (M)}) \leq
d(\overline{ \bigoplus\limits_{k\neq t} Q_k(M) } ,
\ell_2^m) \leq f(n_{t-2})/2,$$
where  $m = {\rm dim}\
(\overline{ \bigoplus\limits_{k\neq t} Q_k(M)} )$.\quad\eop

We refer also to Remark 2.11 below.
\bigskip
{\bf 2.5. Remark.} We show that given any non--decreasing unbounded
functions $g: {\Bbb N} \to {\Bbb R}^+$ and $h: {\Bbb N} \to {\Bbb R}^+$,
$h(1) \geq 3 $, the function $f$ and the sequence $(n_k)$ in the
construction of the space $X$ above can be chosen such that both
of the following are satisfied:

\kohta $1^{\circ}$ $n_k \geq  g(k) n_{k-1} $ for all $k \in {\Bbb N}$,

\kohta $2^{\circ}$ $f(n_k) \leq h(n_{k-2}) $ for all $k \in {\Bbb N},
\ k \geq 2$.

Indeed, we set
$q_0 = 4,\ n_0 = 1$ and $f(1) = 3 $. We assume that  $k \in {\Bbb N}$  and that
$q_t,n_t$  are chosen for  $t < k$ and $f(n) $ is
defined for $n \leq n_{k-1}$.
We define  $q'_k, \ \ 2 < q'_k < q_{k-1}$,
such that
$$n_{k-1}^{1/2-1/q'_k} \leq 1 + \varepsilon/3  \eqno(2.4a)$$
and then  $n_k$ such that $n_k > \max \{2, g(k) \} n_{k-1} $
and such that
$$ h(n_k) \geq (k+2) h(n_{k-1}),   \eqno(2.5a)    $$
$$ n_k^{1/2-1/q'_k} \geq
f(n_{k-1})/2. \eqno(2.5b)$$
Then we define ($f(n) = 3$ for $n_0 < n \leq n_1,$ if $k=1$)
$$f(n) = h(n_{k-2})  \eqno(2.5c) $$
for $n_{k-1} < n \leq n_k $
and, finally, $q_k,\ 2 < q_k \leq q'_k$  such that
$$n_k^{1/2-1/q_k} = f(n_{k-1})/2. \eqno(2.6a)$$
Note that (2.5a,c) imply
$$ f(n_k) \geq k f(n_{k-1})     \eqno(2.7a)  $$
for all $k \geq 3$.

So, in view of (2.4a), (2.5b), (2.6a) and (2.7a), none of the
inequalities (2.4)--(2.6) is changed because of
these extra requirements (except (2.5) for $k\leq 2,$ but this
does not matter), and
hence Proposition 2.4 holds also for the function $f $ and the space $X$
satisfying $1^{\circ}$ and $2^{\circ}$.
\bigskip
Our main application of local distance function is a result which says
that a compact operator  $T \in L(X)$  factors in a ''strong'' sense
through  $Y$  if  $(K,C){\rm -}ld(X^*,Y^* )$  is small enough.  We begin
with a lemma which is in principle known; compare [LT], Proposition
II.5.10.
\bigskip
{\bf 2.6. Lemma.}  {\sl Let  $X$  and  $Y$  be reflexive Banach spaces
such that  $(K,C){\rm -}ld(X^*,Y^*)$  is a finite function for some  $C \geq 1$  and
some  $K: {\Bbb N} \to {\Bbb N}$.  Let $C_1 > C\geq 1$.
Given $C_1$--complemented subspaces
$M \subset X,\ M_Y \subset Y,\ 0 \leq {\rm dim}\ (M) = {\rm dim}\ (M_Y)
=: n < \infty$, projections  $P$  (resp.  $P_Y$)  from  $X$  onto  $M,\
\Vert P\Vert \leq C_1$  (resp.\ from  $Y$  onto  $M_Y,\ \Vert P_Y\Vert
\leq C_1)$  and $m$--codimensional, $m \in {\Bbb N}$, closed subspaces
$X_0 \subset X$  and  $Y_0 \subset Y$, there exist $C$--complemented
subspaces  $N \subset X,\ N_Y \subset Y$  satisfying the following}

\kohta $1^{\circ}$ ${\rm dim}\ (N) = {\rm dim}\ (N_Y) \leq K(n + m)$,

\kohta $2^{\circ}$ $M \subset N,\ M_Y \subset N_Y$,

\kohta $3^{\circ}$ {\sl there exist projections  $Q$  from  $X$  onto  $N$  and  $Q_Y$
from  $Y$  onto  $N_Y$  such that}
\kohta { } $PQ = P,\ P_YQ_Y = P_Y$,
$$({\rm id}_X - Q)(X) \subset X_0,\ ({\rm id}_Y - Q_Y)(Y) \subset Y_0,$$
\kohta { } {\sl and}  $\Vert Q\Vert \leq C_1(C + 2),\
\Vert Q_Y\Vert \leq C_1(C + 2)$

\kohta $4^{\circ}$  $d(N,N_Y) \leq C_1^2(C + 2)^2(K,C){\rm -}ld(X^*,Y^*)(n + m)$.
\bigskip
{\bf Proof.}  Let  $\tilde M := P^*(X^*) \subset X^*,\
\tilde M_Y := P^*_Y(Y^*)
\subset Y^*,\ \hat M := X_0 ^{\perp} \subset X^*$  and  $\hat M_Y :=
Y_0^{\perp} \subset Y^*$.  We have ${\rm dim}\ (\tilde M) =
{\rm dim}\ (\tilde M_Y) = n,\
{\rm dim}\ (\hat M) = {\rm dim}\ (\hat M_Y) = m$.  By the assumption on
the local distance function of  $X^*$  and  $Y^*$  we find
$C$--complemented subspaces  $\tilde N \subset X^*$  and  $\tilde N_Y
\subset Y^*$  such that
$$\tilde M + \hat M \subset \tilde N,\ \tilde M_Y + \hat M_Y \subset
\tilde N_Y, \eqno(2.14)$$
${\rm dim}\ (\tilde N) = {\rm dim}\ (\tilde N_Y) \leq K(n + m)$  and
$d(\tilde N,\tilde N_Y) \leq (1 + \varepsilon)
(K,C){\rm -}ld(X^*,Y^*)(n + m)$, where $\varepsilon > 0 $ is such that
$(1+ \varepsilon ) (CC_1 +C_1 +C)^2 \leq C_1^2 (C+2)^2 $.  Let  $\tilde Q$
and  $\tilde Q_Y$  be projections from  $X^*$  onto  $\tilde N$  and
$Y^*$  onto  $\tilde N_Y$, respectively with norm not greater than  $C$.
We define
$$\hat Q := \tilde Q + P^* - P^*\tilde Q,\ \hat Q_Y := \tilde Q_Y +
P^*_Y - P^*_Y\tilde Q_Y. \eqno(2.15)$$
Then it is elementary to see, using (2.14) that, $\hat Q$  is a
projection from $X^*$  onto  $\tilde N$  which commutes with  $P^*$
and satisfies  $\Vert \hat Q\Vert \leq
CC_1 +C_1 +C\leq C_1 (C + 2)$.  The same holds for
$\hat Q_Y$  with respect to  $Y^*,\tilde N_Y$  and  $P^*_Y$.

Finally, we define using reflexivity the projections  $Q :=
\hat Q^*,\ Q_Y := \hat Q^*_Y$  and the subspaces  $N := Q(X)$  and
$N_Y := Q_Y(Y)$.  To see that  $M \subset N$  let  $x \in M$  and
$y \in X^*$.  Then, by (2.14) and the commutativity of  $P^*$  and
$Q^*$,
$$\eqalign{&\langle Qx,y\rangle = \langle QPx,y\rangle =
\langle x,P^*Q^*y\rangle\cr
\noalign{\vskip 4pt}
&= \langle x,Q^*P^*y\rangle = \langle x,P^*y\rangle = \langle Px,y\rangle\cr
\noalign{\vskip 4pt}
&= \langle x,y\rangle.\cr}$$
Hence, $Qx = x$  for all  $x \in M$, which means that  $M \subset N$.
Concerning the property $3^{\circ}$, the relation  $PQ = P$  follows from
definitions.  Let us prove that  $({\rm id}_X - Q)(X) \subset X_0$.
If  $x \in X$  and  $y \in X_0^{\perp} = \hat M \subset \tilde N \subset
X^*$, then
$$\langle ({\rm id}_X - Q)x,y\rangle = \langle x,y - \hat Qy\rangle
= \langle x,y - y\rangle = 0,$$
since  $\hat Q$  is a projection onto  $\tilde N$.  This implies the
statement.  The norm estimate in $3^{\circ}$  follows from (2.15).

Clearly, the statements for  $N_Y$  etc.\ are proved in the same way.

The distance estimate $4^{\circ}$  follows from  $d(N^*,\tilde N)
\leq CC_1 +C_1 +C , $ $d(N^*_Y,\tilde N_Y)
\leq CC_1 +C_1 +C $  and  $d(\tilde N,
\tilde N_Y) \leq (1 +\varepsilon )
(K,C){\rm -}ld(X^*,Y^*)(n + m)$ and from the choice of
$\varepsilon$.\quad\eop
\bigskip
The main difficulty in applying our method is the following.  We have
finite dimensional subspaces  $M,M_1,\ M_1 \subset M$, and  $N,N_1,\
N_1 \subset N$, of the Banach spaces  $X$  and  $Y$, respectively, such
that  $d(M,N)$  and  $d(M_1,N_1)$  are quite small and such that all the
subspaces are quite well complemented.  However, we do not in general
know if $d(M_2,N_2)$  is small enough
for all "good" complements  $M_2$  and  $N_2$ of $M_1$ and $N_1$ in  $M$
and  $N$, respectively. So,
we give our theorem in two versions.  In the first one we make
assumptions on the spaces  $X$  and  $Y$  which enable us to avoid the
difficulty mentioned above. We will assume that given a Banach space $X$
the following property is satisfied for a suitable non--decreasing
function
$f_X : {\Bbb N} \to {\Bbb R}^+ $
which will be given later:

{\sl (D) There exist constants $d(X) \geq 0 $ and $D(X) > 0$ such that
for each
$n$--dimensional $c$--complemented subspace $M \subset X$ and
projection $P$ from $X$ onto $M$ with $ || P || \leq c $ there
exist an $n$-dimensional
subspace $N$ and a projection $Q$ from $X $ onto $N $
such that }

\kohta $1^{\circ}$ $M \cap N =\{ 0 \}  $,

\kohta $2^{\circ}$ $ || Q || \leq D(X) c^{d(X)}$

\kohta $3^{\circ}$ $ PQ=QP=0$

\kohta $4^{\circ}$ $d(M,N) \leq f_X (n).$

It is not difficult to see  that for example the reflexive
separable ${\cal L}_p$--spaces
satisfy this property with a constant function $f_X$ and $d(X)=2$.
Namely, assume that $X$ is a separable ${\cal L}_{p,\lambda}$--space,
$1 < p< \infty$, $\lambda \geq 1 , $ and that $M \subset X$ is
finite dimensional and that $P $ is a projection from $X$
onto $M$. There exists a subspace $M_1 \supset M$ such that
$d(M_1, {\ell}_p^m ) \leq \lambda $, where $m = {\rm dim } (M_1)$;
let $\phi : M_1 \to {\ell}_p^m$ be an isomorphism with
$ || \phi ||  || \phi^{-1} || \leq \lambda$. It is known that $X$ is
$\lambda$--isomorphic to a subspace of $L_p(0,1)$ (see [LP], Corollary
7.2). On the other hand, every infinite dimensional subspace of
$L_p(0,1)$ which is not isomorphic to ${\ell}_2 $  contains, say, a
$3$--complemented subspace $2$--isomorphic to ${\ell}_p $
(see [KP], the proof of Theorem 2 and Theorem 3, and [LT1], Proposition
1.a.9). So, in view of these remarks ker$(P)$ contains a subspace
$Y$ which is $3\lambda $--complemented in $X$ and $2\lambda
$--isomorphic to ${\ell}_p $    . Let $\psi : Y \to {\ell}_p $ be
an isomorphism with $ || \psi ||  || \psi^{-1} || \leq  2\lambda .$
Let $I_0$  be the natural embedding of $ {\ell}_p^m$ into
$  {\ell}_p$ and let $R_0 $ be  the canonical projection from
$ {\ell}_p  $ onto $ I_0 (  {\ell}_p^m ) $, $\Vert I_0 \Vert =
\Vert R_0 \Vert = 1.$

Now $N=\psi^{-1} I_0 \phi (M) $ is a subspace  of $Y \subset X $ satisfying
$M \cap N = \{ 0 \}$, $ d(M,N) \leq 2 \lambda^2 ,$ and the projection
$$Q= \psi^{-1} I_0 \phi P \phi^{-1} I_0^{-1 } R_0 \psi R ( {\rm id}_X
-P) ,   $$
where $R $ is a projection from $X $ onto $Y$ with $\Vert R \Vert
\leq 3 \lambda$, satisfies $Q(X) =N $, $ PQ=QP=0$, $ \Vert
Q \Vert \leq 12 \Vert P \Vert^2 \lambda^{3}$.

Moreover, if the space $X$ is such that for
a constant $C \geq 1 $, given a finite codimensional closed subspace
$Y$ of $X$ there exists a $C$--complemented subspace $Z$ which is
$C$--isomorphic to $Z$, then (D) is satisfied:
given $M \subset X$ and the projection $P $ from $X$ onto $N$
we take $Y = ({\rm id}_X -P) (X) $ and apply the isomorphy of $Z$
and $X $ to find $N$ analogously to the case
of ${\cal L}_p$--spaces. Furthermore, in Theorem 2.8 it is often
sufficient that a priori only one of the given Banach--spaces satisfy
(D) in full strength. This fits well to the philosophy that
Theorem 2.8 is applicable to "small perturbations" of
"regular" Banach--spaces; see Remark 2.11 below.

\bigskip
{\bf 2.8. Theorem.} {\sl Let  $X$  and  $Y$  be separable, reflexive
Banach spaces such that  $(K,C){\rm -}ld(X^*,Y^*)$  is finite for some
$C \geq 1$  and  $K: {\Bbb N} \to {\Bbb N}$, and let  $T \in L(X)$  be
compact. Choose the sequence $(m_k)_{k=0}^{\infty}$
such that $m_0 = 2,\ m_k \geq  K(4m_{k-1}+2 ) + 1  $
for  $k \in {\Bbb N}$ .
Assume that for some constants $0 \leq \alpha < 1/4 ,
0 \leq \beta < 1/4
$ and  $C' > 0$   the spaces $X$ and $Y $ have property (D) with
$f_X, \ f_Y$ such that for all $k \in {\Bbb N}$
$$f_X(n) \leq C' (a_{m_{k-1}}(T))^{- \alpha} , {\rm for} \ n\leq m_{k+1}$$
$$f_Y(n) \leq C' (a_{m_{k-1}}(T))^{- \beta} , {\rm for} \ n\leq m_{k+1},
 \eqno (2.16) $$
and assume that for all $k \in {\Bbb N}$ the inequality
$$(K,C){\rm -}ld(X^*,Y^*)(4 m_{k+1}) \leq C' k^{-3/4}
C(k+2)^{-5} a_{m_k} (T)^{(-1 + 4 \max \{ \alpha ,\beta \} )/4 },  \eqno(2.17)$$
where  $C(k) := (4 D(C+3))^{(k+1) d^{k}}  $
and $d = \max \{ d(X) ,d(Y),1 \}$, $D = \max \{ D(X)$,$D(Y) \}$
(see property (D)), holds.
Then  $T = T^{(2)} T^{(1)}$,
where $T^{(1)}  \in   L(X,Y)$, $T^{(2)} \in  L(Y,X)$ and
$T^{(1)} (X) $ is dense in $Y$. Moreover, both $T^{(1)}$ and
$T^{(2)} $ are compact.}

\bigskip
Note that the assumptions on $X$ and $Y$ are symmetric
so that the factorization applies as well to compact operators in $Y$.

{\bf Proof.} It is not a restriction to assume that $a_n (T) \leq 1$
for all $n$.

Let for each  $n \in {\Bbb N},\ n > 1$,
the operator  $T_n \in L(X)$
be such that  ${\rm rank}\ (T_n) < n$  and
$$\Vert T - T_n\Vert \leq 2a_n (T) . \eqno(2.18)$$
Let  $(Y_k)^{\infty}_{k=0}$  be a decreasing sequence of closed, finite
codimensional subspaces of  $Y$  satisfying
$${\rm codim}\ (Y_k) = {\rm codim}\ (\bigcap^k_{t=0}\ {\rm ker}\ (T_{m_t})).$$
Let  $(x_n)^{\infty}_{n=1} \subset X$  and
$(y_n)^{\infty}_{n=1} \subset Y$  be sequences of non zero elements
such that  ${\rm sp}\ (x_n) \subset X$  and  ${\rm sp}\ (y_n) \subset Y$
are dense.

We choose the sequences  $(M_k)^{\infty}_{k=1},\ $
$ ({\tilde M }_k)^{\infty}_{k=1},\ $  $ (N_k)^{\infty}_{k=1},\ $
$ ({\tilde N })^{\infty}_{k=1} $
of finite dimensional subspaces of  $X$  or  $Y$  by an inductive method
as follows. We set   $M_0 = \tilde M_0 = N_0 = \tilde N_0 = \{0\},\
P_0 = \tilde P_0 = Q_0 = \tilde Q_0 = 0$.  Assume that  $k \in {\Bbb N}$
and that  $M_n \subset X,\ \tilde M_n \subset X, \
N_n \subset Y,\ \tilde N_n \subset Y$  and the projections  $P_n$  from  $X$
onto  $M_n$  (respectively, $\tilde P_n: X \to \tilde M_n,\ Q_n:
Y \to N_n,\
\tilde Q_n: Y \to \tilde N_n)$  are chosen for  $0 \leq n < k$  such
that the following holds:

\kohta $1^{\circ}$ $M_{n-1} + \tilde M_{n-1} \subset M_n,\ N_{n-1} +
\tilde N_{n-1} \subset N_n$  for all  $1 \leq n < k$,

\kohta $2^{\circ}$  $\Vert P_n\Vert < C(n), \ \Vert \tilde P_n\Vert <
C(n), \
\Vert Q_n\Vert < C(n) ,\ \Vert \tilde Q_n\Vert < C(n)$, for all
$n < k$,

\kohta $3^{\circ}$  $P_n$  commutes with  $P_{n-1}$ and $\tilde P_{n-1}$,
   and  $Q_n$  commutes with
$Q_{n-1}$ and $ \tilde Q_{n-1}$  for all  $1 \leq n < k$,

\kohta $4^{\circ}$  $({\rm id}_X - P_n)(X) \subset
\bigcap\limits^{n-1}_{t=0}\
{\rm ker}\ (T_{m_t}),\ ({\rm id}_Y - Q_n)(Y) \subset Y_{n-1}$  for all
$1 \leq n < k$,

\kohta $5^{\circ}$  $$ d(M_n,\tilde M_n) \leq C' a_{m_{\max \{ n-2,0 \} }}
(T)^{- \alpha},$$
$$ d(N_n,\tilde N_n) \leq C' a_{m_{\max \{ n-2,0 \} }}(T)^{ - \beta}, $$
$$d(M_n,N_n) \leq 16 C(n)^2 (K,C){\rm -}ld(X^*,Y^*)(4m_{n-1}) \ \ {\rm for} \ 1
\leq n < k,$$

\kohta $6^{\circ}$  $P_n\tilde P_n = \tilde P_nP_n = Q_n\tilde Q_n =
\tilde Q_nQ_n = 0$  for all  $n < k$,

\kohta $7^{\circ}$  $x_n \in M_n , \ \ y_n \in N_n$ for all $ 1 \leq n <
k$,

\kohta $8^{\circ}$ dim$(M_n) = $ dim$(\tilde M_n) = $ dim$(N_n) = $
dim$(\tilde N_n) \leq m_n$ for all $n$.

An application of Lemma 2.6 yields the same for  $n \leq k$; more
specifically, we take  $M = M_{k-1} + \tilde M_{k-1},\ M_Y =
N_{k-1} + \tilde N_{k-1},\ X_0 = \bigcap\limits_{t=0}^{k-1}\ {\rm ker}\
(T_{m_t}) ,\ Y_0 = Y_{k-1}$, and for the projections  $P$  and  $P_Y$  we
take  $P_{k-1} + \tilde P_{k-1}$  and  $Q_{k-1} + \tilde Q_{k-1}$,
respectively.  Then we use Lemma 2.6 and first get the spaces
$\hat M_k := N,\ \hat N_k := N_Y$  and the projections  $\hat P_k = Q$
and  $\hat Q_k = Q_Y$  with the  properties $1^{\circ} -4^{\circ}$ of
the lemma (e.g. $\Vert \hat P_k \Vert < 2 C(k-1) (C+2)$
by $2^{\circ}$ above).   Note that, by Lemma
2.6, $\hat P_k$  and  $P_{k-1} + \tilde P_{k-1}$  commute.  Since
$P_{k-1}\tilde P_{k-1} = \tilde P_{k-1}P_{k-1} = 0$, we see that
$$\eqalignno{P_{k-1}\hat P_k &= P_{k-1}(P_{k-1} + \tilde P_{k-1})\hat P_k
= P_{k-1}\hat P_k(P_{k-1} + \tilde P_{k-1})\cr
\noalign{\vskip 4pt}
&= P_{k-1}(P_{k-1} + \tilde P_{k-1}) = P_{k-1}&(2.19)\cr}$$
which means that also  $P_{k-1}$  and  $\hat P_k$, and, similarly,
$\tilde P_{k-1}$  and  $\hat P_k,Q_{k-1}$  and  $\hat Q_k$, as well as,
$\tilde Q_{k-1}$  and  $\hat Q_k$, commute.

We choose the smallest  $n$  (resp.  $m$) such that  $x_n \notin \hat M_k$
(resp.  $y_m \notin \hat N_k)$
and denote  $x^{(k)} := x_n$  (resp.  $y^{(k)} := y_m$).  We have
$({\rm id}_X - \hat P_k)x^{(k)} \neq 0$  and  $({\rm id}_Y - \hat Q_k)
y^{(k)} \neq 0$.  Let  $R_X$  and  $R_Y$  be projections from
$({\rm id}_X - \hat P_k)(X)$  onto  ${\rm sp}\ (({\rm id}_X - \hat P_k)
x^{(k)})$  and from  $({\rm id}_Y - \hat Q_k)(Y)$  onto
${\rm sp}\ (({\rm id}_Y - \hat Q_k)y^{(k)})$, respectively, with norm
one.  We define
$$P_k = \hat P_k + R_X({\rm id}_X - \hat P_k), \ \ Q_k = \hat Q_k + R_Y
({\rm id}_Y - \hat Q_k)$$
and $M_k= P_k (X)$, $ N_k = Q_k (Y)$.
Then clearly

i) $P_k$  commutes with  $\hat P_k$  and  $Q_k$  commutes with  $\hat
Q_k$,

ii) $\Vert P_k\Vert \leq 2\Vert \hat P_k\Vert +1,\ \Vert Q_k\Vert \leq
2\Vert \hat Q_k\Vert +1 $.

Now $1^{\circ}$   above holds for  $n \leq k$. The property $2^{\circ}$
for $P_k$ and $Q_k$ follows from Lemma 2.6:
$$ \Vert P_k \Vert \leq 2 \Vert \hat P_k \Vert + 1
\leq 4 C(k-1)(C+2) +1 . \eqno(2.20)    $$
To verify
$3^{\circ}$ we have, by (2.19) and i)
$$P_{k-1}P_k = P_{k-1}\hat P_kP_k = P_{k-1}\hat P_k = P_{k-1} = P_k
P_{k-1}$$
The same reasoning proves the other cases, too.
Concerning $4^{\circ}$, we have for  $x \in X$
$$({\rm id}_X - \hat P_k)({\rm id}_X - P_k)x = x - \hat P_kx - P_kx +
\hat P_kP_kx = ({\rm id}_X - P_k)x,$$
so that  $({\rm id}_X - P_k)x \in ({\rm id}_X - \hat P_k)(X) \subset
\bigcap\limits^{k-1}_{t=0} {\rm ker}\ (T_{m_t})$, by the choice of
$\hat P_k$.  Similar proofs in  $Y$  show that $3^{\circ}$ and $4^{\circ}$ hold.

The distance estimate for $d(M_k, N_k)$ in $5^{\circ}$ follows from
Lemma 2.6 and the definitions above; note that
the Hahn--Banach--theorem implies the existence of decompositions
$M_k = \hat M_k \oplus M^{(k)}, \ N_k = \hat N_k \oplus N^{(k)},$
where dim$(M^{(k)})= $ dim$(N^{(k)}) =1$, such that the norms of the
corresponding projections do not exceed $2$. This again implies
$d(M_k, N_k) \leq
16  d(\hat M_k , \hat N_k)$, and $4^{\circ}$ of Lemma 2.6 can be applied
to estimate $d(\hat M_k , \hat N_k) $.

We use the assumption on the property (D) of the
spaces  $X$  and  $Y$  and (2.20) to find spaces  $\tilde M_k$  and  $\tilde N_k$,
and projections  $\tilde P_k$  and  $\tilde Q_k$  satisfying the
properties  $2^{\circ}$, $5^{\circ}$ and $6^{\circ}$. It is clear that also
$7^{\circ}$ is satisfied. The property $8^{\circ}$ follows from
definitions.

We denote $P^+_k = P_k + \tilde P_k$ and $Q^+_k = Q_k + \tilde Q_k $
for all $k$. It is easy to see that
$P^+_{k-1 } P^+_k = P^+_{k} P^+_{k-1} = P^+_{k-1}$, and similarly for
$Q^+_k$, hold.

For all  $k$  we denote by  $\psi_k: M_k \to N_k$  an isomorphism
satisfying
$$\Vert \psi_k\Vert \leq 2d(M_k,N_k)^{1/2} , \ \
\Vert \psi_k^{-1}\Vert \leq d(M_k,N_k)^{1/2}$$
and by  $\alpha_k: M_k \to \tilde M_k$  and  $\beta_k: N_k \to
\tilde N_k$  isomorphisms satisfying
$$\Vert \alpha_k\Vert \leq  2 d(M_k,\tilde M_k)^{1/2}, \ \
\Vert \alpha_k^{-1}\Vert \leq d(M_k,\tilde M_k)^{1/2}, $$
$$ \Vert \beta_k\Vert \leq 2 d(N_k,\tilde N_k)^{1/2}, \ \
\Vert \beta_k^{-1}\Vert
\leq d(N_k,\tilde N_k)^{1/2}.$$
We define for all  $k \in {\Bbb N}$  the isomorphisms
$\Phi_k: M_k \oplus \tilde M_k \to N_k \oplus \tilde N_k$,
$$\Phi_k: x + y \mapsto \beta_k\psi_k x  +  \psi_k \alpha_k^{-1}y,$$
where  $x \in M_k,y \in \tilde M_k$.  We have
$${\Phi_k}^{-1}( x + y) = \alpha_k\psi^{-1}_k x  +  \psi^{-1}_k \beta_k^{-1}y,$$
where  $x \in N_k,y \in \tilde N_k$, and
$$\Vert \Phi_k\Vert \leq  4 C(k) d(M_k,N_k)^{1/2}( d(M_k,
\tilde M_k)^{1/2} +  d(N_k,\tilde N_k)^{1/2} ), $$
$$\Vert \Phi_k^{-1}\Vert \leq 2 C(k) d(M_k,N_k)^{1/2}( d(M_k,
\tilde M_k)^{1/2} +  d(N_k,\tilde N_k)^{1/2} ).$$
Let $k \geq 2.$ Note that   $\Phi_k
(M_{k-1} \oplus \tilde M_{k-1}) \subset \tilde N_k$  and
that there exists a projection  $R_k$
from  $N_k \oplus \tilde N_k$  onto  $\Phi_k(M_{k-1} \oplus \tilde M_{k-1} )$  such that
$\Vert R_k\Vert \leq  16 C(k)^3 d(M_k,N_k)( d(M_k,
\tilde M_k)^{1/2} +  d(N_k,\tilde N_k)^{1/2} )^2 $; we can take
$$R_k = \Phi_k P^+_{k-1}\Phi_k^{-1}.\eqno(2.21)$$
Moreover, for $ k \geq 2 $ there exists an isomorphism
$\gamma_k$  from  $N_{k-1} \oplus \tilde N_{k-1}$  onto
$\Phi_k(M_{k-1} \oplus \tilde M_{k-1})$  such that  both $\Vert \gamma_k\Vert $ and
$\Vert \gamma_k^{-1}\Vert $ are not greater than
$$8 C(k)^2 d(M_k,N_k)^{1/2}  ( d(M_k,\tilde M_k)^{1/2} +
d(N_k,\tilde N_k)^{1/2}) d(M_{k-1},N_{k-1} )^{1/2}$$
$$ \times ( d(M_{k-1},\tilde M_{k-1})^{1/2} +
d(N_{k-1},\tilde N_{k-1})^{1/2}); $$
we can take
$$\gamma_k = \Phi_k \Phi_{k-1}^{-1}.\eqno(2.22)$$

We define
$$T^{(1)}x := \sum_{k=1}^{\infty}\ c(k)(({\rm id}_Y - Q^+_{k-1}) +
\gamma_kQ^+_{k-1})\Phi_k(P^+_k - P^+_{k-1})x,\eqno(2.23)$$
where  $P^+_0 = 0, \ \gamma_1 = 0, \ x \in X$  and $c(1) =1,$
$c(k) = k^{-3/2} \Vert \gamma_k\Vert^{-1}
\Vert \Phi_k\Vert^{-1} C(k)^{-2}$ for $k >1,$ and
$$T^{(2)}x = \sum_{k=1}^{\infty}\ c(k)^{-1} T \Phi_k^{-1}
(({\rm id}_Y - R_k) + \gamma_k^{-1}R_k)(Q^+_k - Q^+_{k-1})x,\eqno(2.24)$$
where  $x \in Y,\ Q^+_0 = 0$ and $R_1 = 0.$
It is a direct consequence of the choice of  $c(k)$  that (2.23)
converges absolutely in  $Y$  for all  $x$, and that  $T^{(1)}$  is a
bounded operator. The convergence of the series determined by the right--
hand side of (2.23) is even absolute in $L(X,Y)$ so that $T^{(1)}$ is
compact.

We prove the same for  $T^{(2)}$. Let $k \geq 2.$ We have for all  $x \in Y$
$$R_kQ^+_{k-1}x = \Phi_kP^+_{k-1}\Phi^{-1}_kQ^+_{k-1}x = 0,\eqno(2.25)$$
since  $\Phi_k^{-1}Q^+_{k-1}x \in \tilde M_k$
(see also $6^{\circ}$ and $3^{\circ}$).
Moreover, for  $x \in N_k \oplus \tilde N_k$,
$$Q^+_{k-1}R_kx = 0,\eqno(2.26)$$
since  $R_k x \in \tilde N_k$.  So, denoting  $y(x) := (({\rm id}_Y - R_k)
+ \gamma_k^{-1}R_k)(Q^+_k - Q^+_{k-1})x$  for  $x \in Y$, we get
$$\eqalignno{&P^+_{k-1}\Phi_k^{-1}y(x) = \Phi_k^{-1}\Phi_kP^+_{k-1}\Phi_k^{-1}y(x)
= \Phi_k^{-1}R_ky(x)\cr
\noalign{\vskip 6pt}
&= \Phi_k^{-1}R_k\gamma_k^{-1}R_k(Q^+_k - Q^+_{k-1})x\cr
\noalign{\vskip 6pt}
&= \Phi_k^{-1}R_kQ^+_{k-1}\gamma_k^{-1}R_k(Q^+_k - Q^+_{k-1})x =
0,&(2.27)\cr}$$
by (2.25), since  $\gamma_k^{-1}z \in N_{k-1} \oplus \tilde N_{k-1}$
for all  $z \in R_k(N_k
\oplus \tilde N_k)$.  Clearly, (2.27) implies
$$\Phi_k^{-1}(({\rm id}_Y - R_k) + \gamma_k^{-1}R_k)(Q^+_k - Q^+_{k-1})(Y)
\subset (P^+_k - P^+_{k-1})(X)\eqno(2.28)$$
for all $k \geq 2.$
On the other hand, $(P^+_k - P^+_{k-1})(X) \subset  ({\rm id}_X -
P_{k-1}) (X) \subset \bigcap\limits_{t=0}^{k-2}\
{\rm ker}\ (T_{m_t})$, by the properties $3^{\circ}$, $4^{\circ}$
and $6^{\circ}$
above.  Hence,
for all  $x \in Y$,
$$z_k (x) := \Phi_k^{-1}(({\rm id}_Y - R_k) + \gamma_k^{-1}R_k)(Q^+_k - Q^+_{k-1})x
\in {\rm ker}\ (T_{m_{k-2}}).\eqno(2.29)$$
So, we have by $5^{\circ}$ for all $k \geq 2$
$$\eqalignno{&c(k)^{-1}\Vert Tz_k(x)\Vert =
c(k)^{-1} \Vert (T - T_{m_{k-2}}) z_k(x) \Vert\cr
\noalign{\vskip 6pt}
&\leq 2 k^{3/2}C(k)^2\Vert T - T_{m_{k-2}}\Vert
\Vert \Phi_k\Vert \Vert \Phi_k^{-1}\Vert
\Vert R_k\Vert \Vert \gamma_k\Vert \Vert \gamma_k^{-1}\Vert
\Vert Q^+_k - Q^+_{k-1}\Vert \Vert x\Vert \cr
\noalign{\vskip 6pt}
&\leq 2^{15} k^{3/2} C(k)^{12} a_{m_{k-2}}(T)d(M_k,N_k)^4
(d(M_k,\tilde M_k)^{1/2} + d(N_k,\tilde N_k)^{1/2} )^8 \Vert x \Vert\cr
\noalign{\vskip 6pt}
&\leq 2^{100} k^{3/2} C'^4C(k)^{20}
a_{m_{k-2}}(T)^{1-4\max \{ \alpha, \beta \} }
( (K,C){\rm -}ld(X^*,Y^*)(4m_{k-1}))^{4}
\Vert x \Vert ; &(2.30)\cr}$$
note that we have normalized $a_n (T) \leq 1$ for all $n.$
The inequality (2.17), combined with (2.30) now implies
$$\Vert c(k)^{-1} Tz_k(x)\Vert \leq 2^{100} k^{-3/2} C'^8 \Vert x\Vert $$
for $k \geq 2.$
The constants here do not depend on $k$, so
we see that (2.30) converges for all  $x$  and defines a bounded operator  $T$.
The same inequality even shows that the convergence is absolute in
$L(Y,X)$. Hence, also $T^{(2)}$ is compact.

We show that  $T = T^{(2)}T^{(1)}$. A direct calculation using
$\gamma_k Q^+_{k-1} \subset ({\rm id}_Y - Q^+_{k-1}) (Y) $  shows that for
all $k,t \in {\Bbb N}$ and for all $x \in X$
$$\eqalignno{&(Q^+_t - Q^+_{t-1})(({\rm id}_Y - Q^+_{k-1}) + \gamma_kQ^+_{k-1})
\Phi_k(P^+_k - P^+_{k-1})x\cr
\noalign{\vskip 6pt}
&= \delta_{kt}(({\rm id}_Y - Q^+_{k-1}) + \gamma_kQ^+_{k-1})\Phi_k
(P^+_k - P^+_{k-1})x.&(2.31)\cr}$$
Hence, for  $x \in X$,
$$\eqalignno{&T^{(2)}T^{(1)}x = \sum_{k=1}^{\infty}\ c(k)^{-1} c(k)T
\Phi_k^{-1}(({\rm id}_Y - R_k)\cr
\noalign{\vskip 6pt}
&+ \gamma_k^{-1}R_k) (({\rm id}_Y - Q^+_{k-1}) + \gamma_kQ^+_{k-1})
\Phi_k(P^+_k - P^+_{k-1})x\cr
\noalign{\vskip 6pt}
&= \sum_{k=1}^{\infty}\ T \Phi_k^{-1}({\rm id}_Y - R_k -
Q^+_{k-1} \cr
\noalign{\vskip 6pt}
&+ R_kQ^+_{k-1} + \gamma_k^{-1}R_k({\rm id}_Y - Q^+_{k-1}) + ({\rm id}_Y - R_k)
\gamma_kQ^+_{k-1}\cr
\noalign{\vskip 6pt}
&+ \gamma_k^{-1}R_k\gamma_kQ^+_{k-1})\Phi_k(P^+_k - P^+_{k-1})x.&(2.32)\cr}$$
In view of (2.25) this is equal
to
$$\eqalignno{&\sum_{k=1}^{\infty}\ T\Phi_k^{-1}({\rm id}_Y - R_k - Q^+_{k-1}
+ \gamma_k^{-1}R_k + \gamma_kQ^+_{k-1}\cr
\noalign{\vskip 6pt}
&- R_k\gamma_kQ^+_{k-1} + \gamma_k^{-1}\gamma_kQ^+_{k-1})\Phi_k(P^+_k - P^+_{k-1})x\cr
\noalign{\vskip 6pt}
&= \sum_{k=1}^{\infty}\ T\Phi_k^{-1}({\rm id}_Y - R_k + \gamma_k^{-1}R_k)
\Phi_k  (P^+_k - P^+_{k-1})x.&(2.33)\cr}$$
But we have
$$R_k\Phi_k(P^+_k - P^+_{k-1})x = \Phi_kP^+_{k-1}\Phi_k^{-1}\Phi_k
(P^+_k - P^+_{k-1})x = \Phi_kP^+_{k-1}(P^+_k - P^+_{k-1})x = 0,$$
so that we finally get, using (2.32) and (2.33)
$$T^{(2)}T^{(1)}x = \sum_{k=1}^{\infty}\ T(P^+_k -
P^+_{k-1})x,\eqno(2.34)$$
where  $P^+_0 = 0$.  Note that, as a consequence of the absolute
convergence of the series (2.23) and (2.24), we have absolutely
convergent series in our calculations (2.32), (2.33) and (2.34).  So,
(2.34) implies  $T^{(2)}T^{(1)}x = Tx$  for  $x$  in arbitrary  $M_k$;
now $\overline {\bigcup\limits_k\ M_k} = X$  implies this for all  $x$.

Finally, we show that  $T^{(1)}(X) \subset Y$  is dense.  To this end
it is enough to show that  $N_k \subset T^{(1)}(X)$; see the choice of
$N_k$  and  $(y_k)^{\infty}_{k=1}$  above.  So, assume that  $y \in
(Q^+_k - Q^+_{k-1})(Y)$.  We define
$$x := \Phi_k^{-1}(({\rm id}_Y - R_k)y + \gamma_k^{-1}R_ky).$$
We have clearly  $x \in P^+_k(X)$, and, moreover,
$$P^+_{k-1}x = \Phi_k^{-1}(\Phi_kP^+_{k-1}\Phi_k^{-1}y - \Phi_kP^+_{k-1}
\Phi_k^{-1}R_ky ) + P^+_{k-1}\Phi_k^{-1} \gamma_k^{-1}
R_ky= 0,\eqno(2.35)$$
since  $\Phi_kP^+_{k-1}\Phi_k^{-1} = R_k$ and since
$$\Phi_k^{-1}\gamma_k^{-1}R_k(y) \in \tilde M_k.$$
So, $c(k)^{-1}x \in (P^+_k - P^+_{k-1})(X)$, and we get, by (2.23),
(2.26), by an analogue of (2.31)
and by the assumption  $Q^+_{k-1}y = 0$,
$$\eqalign{&T^{(1)}c(k)^{-1}x = (({\rm id}_Y - Q^+_{k-1}) + \gamma_kQ^+_{k-1})
(({\rm id}_Y - R_k)y + \gamma_k^{-1}R_ky)\cr
\noalign{\vskip 6pt}
&= y - Q^+_{k-1}y - R_ky + Q^+_{k-1}R_ky + \gamma_kQ^+_{k-1}y\cr
\noalign{\vskip 6pt}
&- \gamma_kQ^+_{k-1}R_ky + \gamma_k^{-1}R_ky -
Q^+_{k-1}\gamma_k^{-1}R_ky +
\gamma_k Q^+_{k-1} \gamma_k^{-1}R_ky\cr
\noalign{\vskip 6pt}
&= y - R_ky + \gamma_k^{-1}R_ky - \gamma_k^{-1}R_ky + R_ky = y.\cr}$$
This completes the proof for the density of  $T^{(1)}(X)$  and hence,
for  Theorem 2.8. \quad\eop
\bigskip
The formulation of Theorem 2.8 looks quite technical, so it is
useful to consider some simple numerical examples.

$1^{\circ}$. First note that to get an increasing sequence
$$ \left( k^{-3/4} C(k+2)^{-5} a_{m_k}
(T)^{(-1 + 4 \max \{ \alpha, \beta \} )/4} \right)_{k=1}^{\infty}$$ in (2.17)
it is always necessary to choose the sequence $(m_k)$
such that, roughly speaking, at least
$$ ( a_{m_{k+1}}(T))^{-1} > (a_{m_k}(T) )^{-d}$$
holds for large $k$. To be more exact, if for example
$ ( a_{m_{k+1}}(T))^{-1}   \leq (a_{m_k} (T) )^{-d} $ for all $k$, we have
$$C(k+2)^{-1} a_{m_k} (T)^{-1} \leq 4^{-(k+1) d^k}
a_{m_0} (T)^{-d^k} \to 0,$$
as $k \to {\infty}$, and in this case (2.17) cannot be satisfied.

$2^{\circ}$. Let  $T \in L(X)$  be compact and  $a_n(T) =
n^{-c}$  for some  $c > 0 $.
Now  (2.17) holds, provided for
some  $\varepsilon > 0$ and $\gamma > 0$
$$(K,C){\rm -}ld(X^*,Y^*)(n) \leq (\log\ (n))^{c {(1 - 4 \max \{ \alpha, \beta \}
)/4}  -\varepsilon}$$
and $K(n) \leq e^{\gamma n}$
for all  $n$. (We choose $m_{k+1} = e^{\gamma (4 m_k +2)} +1 .$)

$3^{\circ}$. If $d=1,$ $a_n(T) = (\log\ n )^{-1}$, then (2.17) holds provided
$$(K,C){\rm -}ld(X^*,Y^*)(n) \leq (\log\ (n))^{ (1 - 4 \max \{ \alpha, \beta \}) /8 }$$
for some $c > 0$ and $K(n) \leq n^c$; we choose $m_{k+1} = (7
m_k)^{c'}$,  where
$$ c'= (c+4D( C+3))^{40 /(1 - 4 \max \{ \alpha, \beta \}) }. $$
We get
$$ (K,C){\rm -}ld(X^*,Y^*) (4m_{k+1 }) \leq (\log
(4 m_{k+1}))^{(1 -4 \max \{ \alpha, \beta \}) /8} $$
$$\leq (28 c' \log (m_k)
)^{ (1 - 4 \max \{ \alpha, \beta \}) /8} $$
so that
$$ ( (K,C){\rm -}ld(X^*,Y^*) (4m_{k+1 }) )^{-1}
a_{m_k}(T)^{ (-1 + 4 \max \{ \alpha, \beta \})/4} \geq c_0 ( \log (m_k)
)^{ (1 - 4 \max \{ \alpha, \beta \})/8} \eqno(2.36) $$
for a constant $c_0 > 0$. Using $m_k > m_{k-1}^{c'}$ recursively we get
$m_k > m_0^{c'^k}$, and combining this with (2.36) implies
$$ 2 ( (K,C){\rm -}ld(X^*,Y^*) (4m_{k+1 }) )^{-1}
a_{m_k}(T)^{ (-1 + 4 \max \{ \alpha, \beta \})/4}$$
$$\geq c_0
c'^{k (1 - 4 \max \{ \alpha, \beta \})/8  } \geq c_0 ( c +
4D(C+3))^{5k}. $$
Hence, in view of  the definition of $C(k)$ and $d=1$ we see that
(2.17) holds.
\bigskip
Clearly, Theorem 2.8 is in some cases far from being a necessary
condition for factorization. For example in the case of
a diagonal operator $T$ on ${\ell}_p$, \ \ $Te_n = \lambda_n e_n$,
where $(e_n)_{n=1}^{\infty}$ is the canonical basis and
$(\lambda_n)_{n=1}^{\infty} $ tends to $0$, the operators
$T_n$ connected with the numbers $a_n (T)$ can be chosen such that
the spaces $T_n ( \ell_p )$ are $1$--complemented ${\ell}_p^n$--spaces, and
consequently, the uniformity function $K$ is not always needed.
On the other hand in the case of general operators we do not know how to cope
without the uniformity function.
\bigskip
{\bf 2.9. Lemma.} {\sl Let $X$ and $Y$ be separable reflexive Banach spaces.
Assume that both of them satisfy property (D) with constant functions
$f_X$ and $f_Y$, and that the uniform local distance of $X^*$ and
$Y^*$ is bounded. Then $X$ and $Y $ satisfy the assumptions of Theorem
2.8 for an arbitrary compact $T \in L(X).$ }
\bigskip
To prove this we only need to choose the sequence $(m_k)$ in Theorem
2.8 so rapidly
increasing that the right hand side of (2.17) is for all $k \in {\Bbb N}$
larger than a constant (=the uniform local distance of $X^*$ and
$Y^*$)
\bigskip
{\bf 2.10. Theorem.} {\sl Let $X$ and $Y$ be separable ${\cal L}_p
$--spaces, $1< p < {\infty}$. Then an arbitrary compact $T \in L(X)$
factors through $Y$ as $T= T^{(2)} T^{(1)}$ such that the
image of $X$ in $Y$ is dense and such that both $T^{(1)}$ and $ T^{(2)}$
are compact.}
\bigskip
{\bf Proof.} We refer to Lemma 2.9 and the considerations
after the definition of
property (D) and Remark 2.2,$3^{\circ}$. \quad \eop
\bigskip
{\bf 2.11. Remark.}  We show that if $T \in L(\ell_2)$ is compact,
then
$\ell_2$  and the dual $X^*$ of a
suitably chosen space  $X$  of Proposition 2.4 (with $ \varepsilon = 1/2$)
satisfy the assumptions of Theorem 2.8 with $d(X^*) = 1,$
$D(X^*) = 4,$ $\alpha = 0$, $\beta = 1/8$ and $ m_{k+1} =
K( 4 m_k +2 ) +1.$ The consideration below is just to show that
Theorem 2.8 is applicable to various situations. The conclusion,
Corollary 2.12 itself has also a simpler proof, see Corollary 3.10.
However, see also the remark after Corollary 2.12.

We first show that
there exists a function $ f: {\Bbb N} \to {\Bbb N} $ and a choice
of numbers $n_k$ in the construction of the space $X$ of
Proposition 2.4 such that $(K,C){\rm -}ld(\ell_2,X^{**}) = (K,C){\rm -}ld(\ell_2,X)$
satisfies (2.17). Our aim
is to use Remark 2.5.  We first
choose the function
$g: {\Bbb N} \to {\Bbb N}$ such that $g(k) \geq g(k-1)$ for all
$k >1$, $g(1) \geq 18$ and such that the sequence
$$(k^{-3/4 } C(k+2)^{-5 } a_{g(k)}(T)^{-1/8} )_{k=1}^{\infty},  $$
where $C(k)$ is as in Theorem 2.8 (with $d=1$, $D=4 $ and $C=3/2$), is
increasing and unbounded and the first element is larger than $3$.
We define $h(n) = 3$ for $n \leq g(1)$ and
$$ h(n) = k^{-3/4} C(k+2)^{-5} a_{g(k)} (T)^{-1/8} \eqno(2.37)$$
for $ g(k) \leq n < g(k+1)$.
Then both $h $ and $g$ are non--decreasing unbounded functions and by
Remark 2.5 we can choose the sequence $(n_k)$ and the function $f$
such that $1^{\circ}$ and $2^{\circ}$ of Remark 2.5 hold
and construct the space $X$ satisfying Proposition 2.4 with $ \varepsilon = 1/2.$
By (2.13) we have
$$ n_{k} \leq K(n) \leq 3 n_k , \eqno(2.38) $$
if $n_{k-1} \leq n \leq n_k$.
This implies that
$$ n_k \leq m_k \leq 3 n_k +1. \eqno(2.39)$$
Namely, if for some $k \in {\Bbb N}$
$$n_{k-1} \leq m_{k-1} \leq 3n_{k-1}+1 , \eqno(2.40) $$
we have, by the choice of $g(k) $ above and by the choice $n_k$, Remark
2.5, $1^{\circ}$,
$$ 4 m_{k-1} + 2 \leq 12 n_{k-1} + 6 \leq n_k.$$
Hence, by (2.40) and (2.38) we have
$ n_k \leq  K(4m_{k-1} + 2)   \leq 3 n_k.  $
Since $K(4 m_{k-1 }+2 ) +1 = m_k,$ we get (2.40) for $k$,
and , by induction, this holds for all $k$. (We have
$n_0 = 1, $ $ m_0 = 2$.)
Hence, for all $k \in {\Bbb N}$, by Proposition 2.4 and by $2^{\circ},$
Remark 2.5,
$$ (K,C){\rm -}ld (\ell_2,X) (4 m_{k+1}) \leq f(4 m_{k+1}) \leq f(12 n_{k+1}
+4)
\leq f(n_{k+2} ) \leq h(n_k ) $$
$$\leq h(g(t+1) -1) = t^{-3/4} C(t+2)^{-5} a_{g(t)} (T)^{-1/8}
\leq k^{-3/4} C(k+2)^{-5} a_{n_k} (T)^{-1/8} $$
$$\leq k^{-3/4} C(k+2)^{-5} a_{m_k} (T)^{-1/8},$$
where $t$ is such that $g(t) \leq  n_k < g(t+1)$ (note that
$n_k \geq g(k)$ by Remark 2.5, $1^{\circ}$ so that $t \geq k$ holds).
Hence, (2.17) holds.

Let us then consider property (D). For $\ell_2$ there is nothing
to prove: it satisfies property (D) with $f_{\ell_2} = d({\ell_2})
= 1 $, $D({\ell_2}) \leq 1$, so that we can take $\alpha = 0.$ We note
that
$$ X^* = ( \bigoplus_{k=1}^{\infty} Y^*_k )_{\ell_2}, $$
where $Y_k$ is as in the proof of Proposition 2.4.
Since $Y_k$ is for all $k$ isometric to a subspace of $L_{q_k}$
we see that $Y^*_k$ is a quotient of $L_{p_k}$, where $1/p_k
+ 1/q_k = 1$. So, again by [L], Corollary 5, if $Z$ is an
$n$--dimensional subspace of $Y^*_k$, then
$$  d(Z, {\ell_2}^{{\rm dim}(N)} ) \leq n^{1/p_k - 1/2 }  =
n^{1/2 - 1/q_k}  \eqno(2.41) $$
and there exist a projection from $Y^*_k$, and hence a projection
$R$ from $X^*$, onto $Z$ with
$$\Vert R \Vert  \leq n^{1/2 - 1/q_k }.  \eqno(2.42)$$
Let  $M \subset X^*$  be such that  ${\rm dim}\ (M) \leq m_k \leq
3n_k +  1$, $k> 2 .$ By (2.41) and considerations simpler than in (2.9)--(2.12)
we see that $d(M,{\ell_2}^{{\rm dim}(M)} ) \leq 3 n_k^{1/2 - 1/q_k } +1$.
By (2.6), $2^{\circ}$ of Remark 2.5 and (2.36) we thus get
$$d(M, {\ell_2}^{{\rm dim}(M)} ) \leq  2 f(n_{k-1})  \leq
2 h(n_{k-2})  \leq 2 h(g(t+1)-1) $$
$$\leq  a_{g(t)}(T)^{-1/8} \leq  a_{n_{k-2}}(T)^{-1/8}
\leq  a_{m_{k-2}}(T)^{-1/8}, \eqno(2.43) $$
where $t$ is such that $g(t) \leq  n_{k-2} < g(t+1)$.

On the other hand, let $P$ be a projection from $X^*$ onto $M$. The
subspace $N_0 = {\rm ker}(P) \cap Y^*_{k+2} $ has dimension not
smaller than $ n_{k+2} -m_k$. Let $N$ be any dim$(M)$--dimensional subspace
of $N_0$. We have $d(N_0, \ell_2^{{\rm dim}(N_0)} )
\leq ( 3 n_k+1)^{1/2 - 1/q_{k+2} }
\leq 1 +\varepsilon =3/2,$ see (2.41) and (2.4).
By (2.41), (2.42) and (2.4) we find a projection $R$
from $X^*$ onto $N$ with $\Vert R \Vert \leq  2$. Now  (2.16) follows
from (2.43), and $ Q = R( {\rm id}_X -P )$ is a projection satisfying
the conditions of property (D) with $d(X^*) = 1, \ D(X^*) = 4.$

In fact, we have proved
\bigskip
{\bf 2.12. Corollary.} {\sl Every compact operator $T$ on the
real Hilbert space $\ell_2$ factors
through a separable reflexive Banach space $Y \ \ (=X^* $ above)
without basis such that $T= T_1 T_2,$ where $T_1 \in L(\ell_2,Y),$
$T_2 \in L(Y, \ell_2)$ and $T_1(\ell_2) \subset Y$ is dense.}
\bigskip
{\bf Remark.} The preceding consideration also yields a method to
factorize compact operators on some $\ell_2$-sums of finite dimensional
Banach spaces through the separable Hilbert space. Note that {\sl our
factorization result remains valid, if we replace the spaces
$Y_k$ of (2.7) by arbitrary $n_k$--dimensional subspaces of $L_{q_k}$.}
It is also clear that Theorem 3.2 or similar methods do not work in this
case since all closed subspaces of the Hilbert space are Hilbert spaces.
\bigskip
In the second version of our main result we simply give a condition
for factorization in terms of the solution of problem mentioned after
Lemma 2.6.  To be
more exact, we give the following
\bigskip
{\bf 2.13. Definition.}  Let  $X$  and  $Y$  be Banach spaces and
let $f: {\Bbb N} \to {\Bbb R}^+$, $h: {\Bbb N} \to {\Bbb R}^+$ and
$g: {\Bbb N} \to {\Bbb N}, \ g(n) \leq  n $ for all $n$,
be positive, non--decreasing.  We define (if possible for these
$X,$ $Y,$ $f,$ $g,$ and $h$) the
function  $\Phi(f,g,h): {\Bbb N} \to {\Bbb R}$ as follows.  If  $n \in
{\Bbb N}$, then  $\Phi (f,g,h)(n)$  is the supremum over all pairs
$(M,N)$  of at most $n$--dimensional subspaces  $M \subset X,\ N \subset
Y,$ dim$(M) = $ dim$(N),$
$ d(M,N) \leq f(n)$, and subspaces $M_1 \subset M,\
N_1 \subset N$ satisfying dim$(M_1)= $ dim$(N_1)  \leq g(n)$  and
$d(M_1,N_1) \leq f(g(n))$, of the number
$$\sup \ \{ d(M_2,N_2) \},$$
where the supremum is taken over all complements  $M_2,\ M_2 \oplus M_1 =
M$, and  $N_2,\ N_2 \oplus N_1 = N$ such that the projections
corresponding to these direct sums have norms smaller than $h(g(n))$ .
\bigskip
In general, sufficiently good
estimates for  $\Phi$ are unknown.  The philosophy is
that, if for our spaces for example
$$\Phi (f,g,h) (n) \leq cf(n)^{\alpha} \eqno(2.44) $$
holds for some constants  $c , \alpha >0$ and for some relatively
slowly increasing functions $g$ and $h$, we get reasonable results.

The motivation of Definition 2.13 is the fact that knowing a good
estimate for $\Phi$ makes the proof of our main result quite
straightforward; compare the proofs of Theorems 2.8 and
2.14.
\bigskip
{\bf 2.14. Theorem. } {\sl Let  $X$  and  $Y$ be separable, reflexive
Banach spaces such that  $(K,C){\rm -}ld(X^*$, $Y^*)$  is finite for some
$C \geq 1$  and  $K: {\Bbb N} \to {\Bbb N}$, and let  $T \in L(X)$  be
compact. Let $m_0 = 2,\ m_k \geq  K(3m_{k-1}+2 ) + 1 $
for all $k \in {\Bbb N}$. Let the functions $f: {\Bbb N} \to {\Bbb R}^+$
and $h: {\Bbb N} \to {\Bbb R}^+$ be defined by $f(n) = h(n) = 1$ for
$n\leq 2,$ and
$$ f(n) = 16 C(k(n) )^2 (K,C){\rm -}ld(X^*,Y^*)(3m_{k(n)-1} ), \eqno(2.45)$$
$$ h(n) = C(k(n)) + 1 , \ \ \ {\rm for} \ n > 2 \eqno(2.46) $$
where $k(n) \in {\Bbb N}$ is such that $m_{k(n)-1} < n \leq m_{k(n)}  $
and  $C(k):= (2(C+3))^{k+1},$
and let $g: {\Bbb N} \to {\Bbb N } $ be the largest non--decreasing function
satisfying $g(k) \leq k$ and $ g(m_{k+1} ) = m_k$ for all $k \in {\Bbb N }$.}

{\sl Assume that for a constant $C'$ and
for all  $k \in {\Bbb N}$ the inequality
$$ \Phi (f,g,h)(n) \leq (k+1)^{-3} C' C (k+1)^{-1} a_{m_{k-1}}^{-1} (T)
, \ {\rm for \ all } \ n \leq m_{k+1}, \eqno(2.47)$$
 holds ($\Phi$ with respect to the spaces $X$ and $Y$). Then
$T = T^{(2)} T^{(1)}$,
where $T^{(1)}  \in   L(X,Y)$, $T^{(2)} \in  L(Y,X)$ and
$T^{(1)} (X) $ is dense in $Y$. Moreover, the operators
$T^{(1)}$ and $T^{(2)}$ are compact.}
\bigskip
Note that if the estimate (2.44) holds, then the condition
$$ (K,C){\rm -}ld(X^*,Y^*) (3m_k) \leq (k+1)^{-3/\alpha}
C(k+1)^{-2-1/ \alpha}a_{m_{k-1}}(T)^{-1 / \alpha } \eqno(2.48) $$
implies (2.47), and hence, the operator $T$ factors. (Combine (2.44),
(2.45) and (2.48).) On the other hand, (2.48) resembles (2.17) so
we get numerical examples like those after the proof of Theorem 2.8, but
with a bit better $C(k).$ Of course, applying Theorem 2.14 one does not
need to take care of property (D).
\bigskip
{\bf Proof.} For each  $n \in {\Bbb N},\ n > 1$, we choose an
operator $T_n \in L(X)$  satisfying  ${\rm rank}\ (T_n) < n$  and
$$\Vert T - T_n\Vert \leq 2a_n(T).\eqno(2.49)$$
Let  $(Y_n)^{\infty}_{n=0}$  be a sequence of closed subspaces of  $Y$
such that   $Y_n
\supset Y_{n+1}$  and  ${\rm codim}\ (Y_n) = {\rm codim}\
(\bigcap\limits_{k=0}^n {\rm ker}\ (T_{m_k}))$.

\indent Let  $(x_n)^{\infty}_{n=1} \subset X$  and
$(y_n)^{\infty}_{n=1} \subset Y$  be sequences of non zero elements
such that  ${\rm sp}\ (x_n) \subset X$  and  ${\rm sp}\ (y_n) \subset Y$
are dense.

We apply Lemma 2.6 inductively as follows.
We set $M_0 = N_0 = \{0\}$, \ \ $P_0 = Q_0 = 0$.
Assume that  $k \geq 1$
and that for $0 \leq n < k$ the projections
$P_n \in L(X)$  and  $Q_n \in L(Y)$  are
defined such that for all $0 < n < k$

\kohta $1^{\circ}$  $\Vert P_n\Vert \leq C(n),\Vert Q_n\Vert \leq C(n)$,

\kohta $2^{\circ}$ $P_n (X) \supset P_{n-1}(X), \ \ Q_n (Y) \supset
Q_{n-1}(Y),$

\kohta $3^{\circ}$ $$({\rm id}_X - P_n )(X) \subset ({\rm id}_X -
P_{n-1}) (X) \cap \bigcap_{t=0}^{n-1} {\rm ker} (T_{m_t}), $$
$$ ({\rm id}_Y - Q_n )(Y) \subset ({\rm id}_Y -
Q_{n-1}) (Y) \cap Y_{n-1}, $$

\kohta $4^{\circ}$ dim $(M_n) = $ dim $(N_n) \leq m_n$,

\kohta $5^{\circ}$  $d(M_n,N_n) \leq 16 C(n)^2 KC\ell
df(X^*,Y^*)(3m_{n-1})$

\kohta $6^{\circ}$ $x_{n} \in M_{n}, \ \ y_{n} \in N_n,$

\noindent where  $M_n := P_n(X)$  and  $N_n := Q_n(Y)$.  Then applying
Lemma 2.6
with  $X,M_{k-1},$ $\bigcap\limits_{t=0}^{k-1} {\rm ker}\ (T_{m_t}),$
$Y,N_{k-1}$  and
$Y_{k-1}$ as $X,M,X_0,Y,M_Y$  and  $Y_0$,
and $P_{k-1}$ and $Q_{k-1}$ as $P$ and $P_Y$
we get projections $Q = \hat P_k $ and $Q_Y = \hat Q_k $
onto the subspaces $\hat M_{k} = N \subset X$ and
$\hat N_{k} = N_Y \subset Y$, respectively, satisfying the properties
mentioned in Lemma 2.6. So, we have $|| \hat P_k || \leq 2^k (C+3)^k
(C+2)$  and  $|| \hat Q_k || \leq 2^k (C+3)^k (C+2).$

We choose the smallest  $n$  (resp.  $m$) such that  $x_n \notin \hat M_k$
(resp.  $y_m \notin \hat N_k)$
and denote  $x^{(k)} := x_n$  (resp.  $y^{(k)} := y_m$).  We have
$({\rm id}_X - \hat P_k)x^{(k)} \neq 0$  and  $({\rm id}_Y - \hat Q_k)
y^{(k)} \neq 0$.  Let  $R_X$  and  $R_Y$  be projections from
$({\rm id}_X - \hat P_k)(X)$  onto  ${\rm sp}\ (({\rm id}_X - \hat P_k)
x^{(k)})$  and from  $({\rm id}_Y - \hat Q_k)(Y)$  onto
${\rm sp}\ (({\rm id}_Y - \hat Q_k)y^{(k)})$, respectively, with norm
one.  We define
$$P_k = \hat P_k + R_X({\rm id}_X - \hat P_k), \ \ Q_k = \hat Q_k + R_Y
({\rm id}_Y - \hat Q_k)$$
and $M_k = P_k(X), $ $N_k =Q_k (Y)$.
Then clearly

i) $P_k$  commutes with  $\hat P_k$  and  $Q_k$  commutes with  $\hat
Q_k$,

ii) $\Vert P_k\Vert \leq 2\Vert \hat P_k\Vert + 1,\ \Vert Q_k\Vert \leq
2\Vert \hat Q_k\Vert + 1$.

Now it is straightforward  to see using Lemma 2.6 that properties
$1^{\circ}$--$6^{\circ}$ are satisfied
for  $k$ instead of $k-1$. We just remark that $ ({\rm id}_X - P_k )(X) \subset
({\rm id}_X - P_{k-1}) (X)  $ follows from the commutativity of $P_k$
and $P_{k-1}$, which again is a consequence of Lemma 2.6 and i).
Note also that
the Hahn--Banach--theorem implies the existence of decompositions
$M_k = \hat M_k \oplus M^{(k)}, \ N_k = \hat N_k \oplus N^{(k)},$
where dim$(M^{(k)})= $ dim$(N^{(k)}) =1$, such that the norms of the
corresponding projections do not exceed $2$. This again implies
$d(M_k, N_k) \leq
16 d(\hat M_k , \hat N_k)$. Lemma 2.6 yields an estimate for
$d(\hat M_k , \hat N_k)  $, so we get $5^{\circ}$.

Using Definition  2.13, (2.45), (2.46) and
the properties $1^{\circ}$ and $4^{\circ}$
above we define for all  $k \in {\Bbb N}, \ k > 1$  the
isomorphism  $\psi_k: (P_k - P_{k-1})(X) \to (Q_k - Q_{k-1})(Y)$  such
that
$$\Vert \psi_k\Vert \Vert \psi_k^{-1}\Vert \leq \Phi (f,g,h) (m_k).
\eqno(2.50)$$
(We take $n= m_k $ and the spaces
$M_k,$ $ M_{k-1},$ $(P_k - P_{k-1}) (M_k)$ correspond to
$M, \ M_1$ and $M_2$ in Definition 2.13; similarly for $N_k$ etc.)
Let $\psi_1 : M_1 \to N_1$ be an isomorphism such that
$ \Vert \psi_1 \Vert \Vert \psi_1^{-1}\Vert $ = $ d(M_1 , N_1). $

We now define
$$T^{(1)}x = \sum_{k=1}^{\infty} k^{-3/2}C(k)^{-1}\Vert \psi_k\Vert^{-1}
\psi_k(P_k - P_{k-1})x,\eqno(2.51)$$
where  $P_0 = 0$  and  $x \in X$,
$$T^{(2)}x = \sum_{k=1}^{\infty} k^{3/2} C(k)\Vert \psi_k\Vert T\psi_k^{-1}
(Q_k - Q_{k-1})x,\eqno(2.52)$$
where  $Q_0 = 0$  and  $x \in Y$.  The facts that  $T^{(1)} \in L(X,Y)$
and that $T^{(1)}$ is compact
are consequences of  $\Vert P_k\Vert \leq C(k),\ k \in {\Bbb N}$.
Concerning  $T^{(2)},$  we have by $2^{\circ}$ for $k \leq 2$ and
for all  $x \in Y$
$$\psi_k^{-1}(Q_k - Q_{k-1})x \in  (P_k - P_{k-1})(X) =
({\rm id}_X -P_{k-1})P_k (X) \subset
{\rm ker}\ (T_{m_{k-2}}).$$
Hence, by (2.49), for $k\geq 2$
$$\eqalign{&\Vert T\psi_k^{-1}(Q_k - Q_{k-1} ) x\Vert = \Vert (T -
T_{m_{k-2}})   \psi_k^{-1}(Q_k - Q_{k-1})x\Vert\cr
\noalign{\vskip 6pt}
&\leq \Vert T - T_{m_{k-2}}\Vert \Vert \psi_k^{-1}\Vert \Vert Q_k - Q_{k-1}\Vert
\Vert x\Vert\cr
\noalign{\vskip 6pt}
&\leq 4C(k)a_{m_{k-2}}(T)\Vert \psi_k^{-1}\Vert \Vert x\Vert.\cr}$$

Now the choice of  $\psi_k$, (2.50) and (2.47) imply for $k \geq 2$
$$\eqalignno{& \left\Vert  C(k) \Vert \psi_k \Vert T\psi_k^{-1}(Q_k -
Q_{k-1})x \right\Vert\cr
\noalign{\vskip 6pt}
&\leq 4C(k)a_{m_{k-2}}(T)\Vert \psi_k\Vert \Vert \psi_k^{-1}\Vert
\Vert x\Vert\cr
\noalign{\vskip 6pt}
&\leq 4C(k) a_{m_{k-2}} (T) \Phi (f,g,h) (m_k)  \cr
\noalign{\vskip 6pt}
&\leq 4 k^{-3} C' \Vert x\Vert. & (2.53)\cr}$$
So (2.53) shows that
(2.52) converges absolutely for all  $x$,   $T^{(2)} \in
L(Y,X)$ and $T^{(2)}$ is even compact.

We have
$$\eqalignno{&T^{(2)}T^{(1)}X = \sum_{k=1}^{\infty} k^{3/2}
C(k)\Vert \psi_k\Vert
T\psi_k^{-1}(Q_{k} - Q_{k-1}) (\sum_{t=1}^{\infty} t^{-3/2} C(t)^{-1}
\Vert \psi_t\Vert^{-1}\psi_t(P_{t} - P_{t-1})(x))\cr
\noalign{\vskip 8pt}
&= \sum_{k=1}^{\infty}\ \sum_{t=1}^{\infty}\  k^{3/2} t^{-3/2}
C(k)C(t)^{-1}
\Vert \psi_k\Vert \Vert \psi_t\Vert^{-1}T\psi_k^{-1}(Q_{k} - Q_{k-1})
\psi_t(P_{t} - P_{t-1})(x). & (2.54) \cr}$$
Here
$$\eqalign{&T\psi_k^{-1}(Q_k - Q_{k-1})\psi_t(P_{t} - P_{t-1})(x)\cr
\noalign{\vskip 6pt}
&= \delta_{kt}T(P_{t} - P_{t-1})(x)\cr}$$
so that (2.24) is equal to
$$\sum_{t=1}^{\infty} T(P_{t} - P_{t-1})(x).$$
So  $T^{(2)}T^{(1)}x = Tx$  for  $x \in M_k$   for
all  $k \in {\Bbb N}$, and hence also for all
$$x \in \overline{\bigcup^{\infty}_{k=1} M_k} = X. $$

The density of $T^{(1)} (X) \subset Y$ follows from the choice of
the sequence $(y_k)$ and the fact that $y_k \in T(M_{k+1})$ for all
$k \in {\Bbb N}$.
\quad\eop
\bigskip
{\bf 2.15. Corollary.} {\sl $1^{\circ}$. Let $T \in L (\ell_2)$ be compact,
$a_n (T) \leq (\log n)^{-1 -\varepsilon}$ for some $0 < \varepsilon <1.$
Then $T $ factors as $T = T_2 T_1$ through every weak Hilbert space
$X$ such that $T_1 ( \ell_2) \subset X$ is dense.

$2^{\circ}$. If  $X$ is a  weak Hilbert space, $T \in L (X)$ is compact
and $a_n (T) \leq (\log n)^{-1 -\varepsilon}$, then $T$ factors
through a Hilbert space.}
\bigskip
Note that to prove this we need to use the unpublished results mentioned
in Remark 2.2.$8^{\circ}$.
\bigskip
{\bf Proof.} Let $X$ be weak Hilbert, let $C=2$ and $K: {\Bbb N}
\to {\Bbb N}$, $K(n) = c'n,$ where $c' \geq 1$ is such that
$(K,C){\rm -}ld(X, \ell_2) (n) \leq c_0 \log (n+1)  $ as in Remark
2.2.$8^{\circ}$. Let $c>0$ be such that if $M \subset X$ is an
$n$--dimensional subspace, then $d(M,\ell_2^{{\rm dim}(M)} ) \leq
c \log (n+1) ,$ see Remark 2.2.$8^{\circ}$. We define $\alpha =
10^{2 / \varepsilon }$ and $m_k = \max \{ e^{\alpha^k  } ,\
7c' m_{k-1} \}$ for all $k \in {\Bbb N};$ for large $k$ we thus
have $m_k =  e^{\alpha^k  }.$

Assume now that $k \in {\Bbb N}$ and $n \leq m_{k+1}$. By the choice
of $c$ and the definion of $\Phi$ we have the following trivial
estimate for $\Phi (f,g,h) (n)$ (with $X=X,$ $ Y= \ell_2$)
which in fact does not depend on $f,$ $g$ and $h:$
$$
\Phi (f,g,h)(n) \leq c \log (m_{k+1} +1).
$$
On the other hand $(a_{m_{k-1}} (T))^{-1} \geq
(\log m_{k-1} )^{1 + \varepsilon}$ so that we get for some constants
$c_i > 0$, $i=1,2,3,$ for  large $k $ and $n \leq m_{k+1}$
$$
(\Phi (f,g,h) (n) )^{-1} a_{m_{k-1}}^{-1}(T) \geq
c^{-1}_1 \alpha^{-k-2} \alpha^{(k-1)(1+ \varepsilon)} $$
$$\geq c_2 \alpha^{k \varepsilon} = c_2 10^{2k} \geq
c_3 10^{k}C(k+1),
$$
see the choice of $C$ and $C(k).$ This shows that (2.47) holds.
So, Theorem 2.14 applies to prove our Corollary in both cases.
\quad\eop
\bigskip
We now turn to a study of locally convex spaces.
Factorization theorems in Banach spaces can often be
used to find new systems of local Banach spaces of locally convex
spaces.  For example, given a Schwartz space  $E$  as a projective
limit of nonreflexive spaces we can use the factorization result of
[DFJP]
which says that a weakly compact operator always factors through a
reflexive space, to find a system of reflexive local Banach spaces on
$E$.  However, not all factorization results are useful in this respect.
To get a trivial counterexample, let us consider a Banach space
$(X,\Vert \cdot \Vert )$ as a Fr\'echet space with the system
$(p_k)_{k=1}^{\infty}$, $p_k = k \Vert \cdot \Vert$, of seminorms.
Given any Banach space $Y$ the linking map (the identity operator on
$X$) between the local Banach spaces $X_{p_{k+1}}$ and $X_{p_{k}}$
factors trivially through $X \times Y$. However, it is not possible
that $X$ can have a system of local Banach spaces isomorphic to
$X \times Y$
unless $X$ is isomorphic to $X \times Y.$

The situation is different, if the factorization is dense.
To see this, assume that  $(E,(p_{\alpha})_{\alpha \in A})$  is a locally
convex space.
Let  $\alpha,\beta \in A,\ p_{\alpha}
\geq p_{\beta}$, and let  $T_{\alpha,\beta}: E_{p_{\alpha}} \to
E_{p_{\beta}}$
be the canonical mapping induced by the identity operator on  $E$.  If
$T_{\alpha \beta}$  factors as  $T_{\alpha \beta} = T^{(2)}T^{(1)}$
through some Banach spaces  $Y$  such that  $T^{(1)}(E_{\alpha}) \subset Y$
is dense, then there exists a continuous seminorm  $q$  on  $E$  such
that  $p_{\alpha} \leq q \leq p_{\beta}$  and  $E_q \cong Y$:
we can take $q(x) = \Vert T^{(1)} \psi x \Vert_{Y} $, where
$\psi$ is the quotient mapping from $E$ onto $E/$ker$(p_{\alpha})$. Of
course, this does not necessarily hold without the assumption on the
density of  $T^{(1)}(E_{p_{\alpha}})$.

Recall that in the case of Schwartz spaces we can find for all  $\alpha \in A$
an index  $\beta$  such that  $T_{\alpha \beta}$  is compact.

Assume that there exists  a separable,
reflexive Banach space  $X$ such that
for all  $\alpha \in A,\ E_{p_{\alpha}} \cong X$.
In view of the preceding remarks it is now clear that, given a Schwartz space as
above, we can use Theorems 2.8 and 2.14 to find new systems of seminorms
$(q_{\beta})_{\beta \in B}$  on  $E$  such that the local Banach spaces
of this systems are not necessarily isomorphic to  $X$.

We could formalize this statement as a corollary but we do not want to
repeat the many assumptions of Theorems 2.8 and 2.14.  We just give some special
cases.
\bigskip
{\bf 2.16. Corollary.} {\sl Let  $E$  be a Schwartz space as above
such that  $X$  is isomorphic to a ${\cal L}_p$--space, $1 < p < \infty$.
Given any ${\cal L}_p$--space  $Y$  the space  $E$  has a system of local
Banach spaces isomorphic to}  $Y$.
\bigskip
To prove this we use Theorem 2.10.
\bigskip
{\bf 2.17. Corollary.}  {\sl Let  $E$  be a hilbertizable
Fr\'echet--Schwartz space over ${\Bbb R}$, i.e., a real
Fr\'echet--Schwartz space having a
system of local Banach spaces isomorphic to  $\ell_2$.  Then there
exists a system
of local Banach spaces $(E_{p_k})_{k=1}^{\infty}$ such that none of the
spaces $E_{p_k}$ has a basis.}
\bigskip
This follows from Corollary 2.12.

I think it should be possible to choose the spaces $E_{p_k} $ isomorphic
to each other.
%
\bigskip
{\bf 3. Remarks on dense factorizations.}  In this section we present
some relatively simple, but efficient, remarks on dense
factorizations of operators in Banach
spaces.  The results can be combined for example with the considerations
in the previous section to get much more new examples of systems of
local Banach spaces in locally convex spaces.

We recall the trivial
fact (see also the remark above Corollary 2.16)
that given a Banach Space  $X$
the identity operator on  $X$  factors through  $X \times Y$  for any
Banach space  $Y$, but that this factorization can in general not be
done such that the image of  $X$  in  $X \times Y$  is dense, even if
$Y$  is separable.  We shall see in this section that the situation
changes dramatically, if we consider for example a compact operator on
$X$.

We begin with a lemma in the spirit of [V], Theorem 1: to prove the
lemma we use the existence of a total, bounded biorthogonal system in a
separable Banach space, as in [V].
\bigskip
{\bf 3.1. Lemma.} {\sl Let  $X$  be a Banach space with a normalized
basis  $(e_n)^{\infty}_{n=1}$, let  $Y$  and  $Z$  be Banach spaces,
$Z$  separable, and let  $T \in L(X,Y)$  satisfy for some}  $\varepsilon
> 0$
$$\Vert Te_n\Vert \leq n^{-2-\varepsilon}\eqno(3.1)$$
{\sl for all  $n$.  Then  $T$  factors as  $T = T^{(2)}T^{(1)}$  through
$Z$  such that  $T^{(1)}(X) \subset Z$  is dense and such that  $T^{(1)}$
and  $T^{(2)}$  are compact.}
\bigskip
{\bf Proof.}  Let  $((x_n)_{n=1}^{\infty},(y_n)^{\infty}_{n=1})$  be a
total, fundamental, biorthogonal system in  $Z$, as in [LT1], Theorem
1.f.4.  Recall that  $((x_n)^{\infty}_{n=1},(y_n)^{\infty}_{n=1})$
satisfies the following properties: $(x_n)^{\infty}_{n=1} \subset Z$
is total, $\langle x_n,y_m\rangle = \delta_{mn}$  for  $n,m \in {\Bbb N}$,
and  $\sup\limits_{n\in {\Bbb N}}\ \Vert x_n\Vert\Vert y_n\Vert \leq 20$,
so that we may assume  $\Vert x_n\Vert \leq 1,\ \Vert y_n\Vert \leq 20$
for all  $n$.

Let  $(e^*_n)^{\infty}_{n=1} \subset X^*$  be the sequence of the
coefficient functionals of  $(e_n)^{\infty}_{n=1}$; we have
$\sup\limits_n\ \Vert e^*_n\Vert < C < \infty$  for a constant
$C > 0$.

We define
$$T^{(1)}x = \sum_{n=1}^{\infty} n^{-1-\varepsilon/2}
\langle x , e^*_n\rangle x_n\eqno(3.2)$$
for  $x \in X$, and
$$T^{(2)}x = \sum_{n=1}^{\infty} n^{1+\varepsilon/2}
\langle x,y_n\rangle Te_n$$
for  $x \in Z$.  The boundedness of the sequences
$(e^*_n)^{\infty}_{n=1}$  and  $(x_n)^{\infty}_{n=1}$  implies
$$\Vert T^{(1)}x\Vert \leq \sum_{n=1}^{\infty} n^{-1-\varepsilon/2}
\Vert x\Vert\Vert e^*_n\Vert\Vert x_n\Vert \leq C'\Vert x\Vert$$
for a constant  $C' > 0$, so that  $T^{(1)}$  is a bounded and even
a compact operator.  Moreover, $T^{(1)}(X)$  is dense in  $Z$  since
$(x_n)_{n=1}^{\infty}$  is total in  $Z$.  Similarly, by (3.1), for
$x \in Z$,
$$\Vert T^{(2)} x \Vert \leq \sum_{n=1}^{\infty} n^{1+\varepsilon /2}
\Vert x\Vert \Vert y_n\Vert \Vert Te_n\Vert \leq \sum_{n=1}^{\infty}
20n^{-1-\varepsilon /2}\Vert x\Vert.$$
Hence, also  $T^{(2)}$  is bounded and compact.

We have for all  $x \in X$
$$\eqalign{T^{(2)}T^{(1)}x &= \sum_{n=1}^{\infty}\ \sum_{m=1}^{\infty}
n^{1+\varepsilon /2}m^{-1-\varepsilon /2}\langle x,e^*_m\rangle
\langle x_m,y_n\rangle Te_n\cr
\noalign{\vskip 6pt}
&= \sum_{n=1}^{\infty} n^{1+\varepsilon /2}n^{-1-\varepsilon /2}
\langle x , e^*_n\rangle Te_n = Tx,\cr}$$
since  $\langle x_m,y_n\rangle = \delta_{m,n}$.\quad\eop
\bigskip
The following result provides a trick to make dense factorizations.
We again refer to [LT1] for the terminology of bases and basic sequences
in Banach spaces.
\bigskip
{\bf 3.2. Theorem.}  {\sl Let  $X$  and  $Y$  be Banach spaces, let
$T \in L(X,Y)$  be compact and assume that  $X$  has a complemented
unconditional basic sequence.  If $Z$  is an arbitrary separable
Banach space, then  $T$  factors as  $T = T^{(2)}T^{(1)}$  through
$X \times Z$  such that  $T^{(1)} (X) \subset X \times Z$  is dense and
$T^{(2)}$  is compact.}
\bigskip
{\bf Proof.}  Let  $(e_n)^{\infty}_{n=1} \subset X$  be a normalized
complemented unconditional basic sequence.  By compactness,
$(Te_n)^{\infty}_{n=1}$  has a convergent subsequence
$(Te_{n_k})^{\infty}_{k=1}$, and, hence, the sequence
$$(Te_{n_{2k}} - Te_{n_{2k+1}})^{\infty}_{k=1}$$
converges to  0  in  $Y$.  Let  $(m_k)^{\infty}_{k=1}$  be a subsequence
of  $(n_k)$  such that for all  $k \in {\Bbb N}$
$$\Vert Te_{m_{2k}} - Te_{m_{2k+1}}\Vert \leq k^{-3},$$
and let us denote  $f_k = e_{m_{2k}} - e_{m_{2k+1}}$  for  $k \in {\Bbb N}$.
Since  $(e_n)^{\infty}_{n=1}$  is an unconditional, complemented basic
sequence we see that  $(f_n)^{\infty}_{n=1}$  is a complemented block
basic sequence of  $(e_n)$  satisfying
$$\Vert Tf_n\Vert \leq n^{-3}\eqno(3.3)$$
for all  $n \in {\Bbb N}$.  Moreover, we have
$$2 \geq \Vert f_k\Vert = \Vert e_{m_{2k}} - e_{m_{2k+1}}\Vert \geq
C^{-1} \Vert e_{m_{2k}}\Vert = C^{-1} $$
for all  $k$, where  $C > 0 $  is the basis constant of
$(e_n)^{\infty}_{n=1}$.

By Lemma 3.1 and (3.3), the restriction of  $T$  to the complemented subspace
$X_1 := {\rm sp}\ (f_n\mid n \in {\Bbb N})$  of  $X$  factors densely
through  $X_1 \times Z$.  Let us denote this factorization by
$T|_{X_1} = R^{(2)}R^{(1)}$, where the operators  $R^{(1)} \in
L(X_1,X_1 \times Z)$  and  $R^{(2)} \in L(X_1 \times Z,X)$  are compact.
We define
$$T^{(1)} = R^{(1)} \oplus {\rm id}_{X_2},\ \ T^{(2)} =
R^{(2)}P + TQ,\eqno(3.4)$$
where  $X_2$  is a closed subspace of  $X$  such that  $X_1 \oplus X_2 = X$,
and  $P$  (respectively, $Q$) is the natural projection from
$(X_1 \oplus X_2) \times Z$  onto  $X_1 \times Z$  (respectively, onto
$X_2$).  Clearly, $T^{(2)}$  is compact as a sum of two compact
operators, see Lemma 3.1.\quad\eop
\bigskip
{\bf Examples.} $1^{\circ}$ The straightforward application of Theorem 3.2 to the Schwartz space
$(E,(p_{\alpha})_{\alpha \in A})$  yields the result that if the spaces
$E_{p_{\alpha}}$  have complemented unconditional basic sequences, then we
can find a system of local Banach spaces  $(E_{q_{\alpha}})_{\alpha \in A}$
where each  $E_{q_{\alpha}}$  is an arbitrary separable Banach space
containing  $E_{p_{\alpha}}$  as a complemented subspace.  This result is
quite natural in view of the facts that each Schwartz space contains a
nuclear subspace and that in nuclear spaces the geometry of local
Banach spaces may be chosen arbitrarily, [V].  However, these known
facts do not imply our results.

$2^{\circ}$ Theorem 3.2 also leads to many examples of the following
type: Let  $(E_n)^{\infty}_{n=1}$  be a family of Banach spaces and let
$\ell_p((E_n)^{\infty}_{n=1})$  be the Banach space of $E_n$--valued,
$\ell_p$--summable sequences, where  $1 \leq p < \infty$.  The space
$\ell_p((E_n))$ contains a complemented copy of  $\ell_p$  (fix one
vector  $0 \neq e_n \in E_n$  for all  $n$  and use the Hahn--Banach
theorem to find for all  $n \in {\Bbb N}$  a projection from  $E_n$
onto  ${\rm sp}\ (e_n))$.  By Theorem 3.2, every compact operator on
$\ell_p$  factors densely through  $\ell_p((E_n))$.  Generalizations
and applications to Schwartz spaces are left to the reader.
\bigskip
{\bf 3.3. Corollary.} {\sl Every compact
operator  $T: \ell_1 \to \ell_1$
factors through an arbitrary separable ${\cal L}_1$--space
$Y$
such that the image of  $\ell_1$  in  $Y$  is dense.}
\bigskip
{\bf Proof.} Follows immediately from Theorem 3.2, since each
${\cal L}_1$--space contains a complemented copy of  $\ell_1$,
[LP], Proposition 7.3.
\quad\eop
\bigskip
Corollary 3.3 complements Theorem 2.10.

In the case of $\ell_2$  a much
stronger result holds.
\bigskip
{\bf 3.4. Theorem. } {\sl Let $T \in L( \ell_2 )$
be compact and let $X$ be an arbitrary separable ${\cal L}_p$--space,
where $1 < p < \infty$.
The operator $T$ factors as $T= T^{(2)} T^{(1)}$ through
$X$ such that $T(\ell_2) \subset X$ is dense and both
$T^{(1)} \in L(\ell_2 , X )$ and $T^{(2)} \in L(X, \ell_2)$ are
compact. }
\bigskip
{\bf Proof.} Using for example the polar decomposition we can find compact
operators $S^{(i)} \in L(\ell_2)$, $i =1,2$, such that $T = S^{(2)}
S^{(1)}$ and such that  each $S^{(i)}$ has a dense range.
The space $L_p(0,1)$
has a complemented subspace $Y $ isomorphic to $\ell_2$ (see, e.g.
[LT], p.215). By Theorem
3.2, both $S^{(1)}$ and $S^{(2)}$ factor through $L_p(0,1)$,
$$
S^{(1)} = S^{(1,2)} S^{(1,1)} \  , \ \ S^{(2)} = S^{(2,2)}
S^{(2,1)},
$$
such that $S^{(i,1)}   \subset L_p(0,1)$ is dense
and $ S^{(i,2)} $ is compact for $ i = 1,2.$ By Theorem 2.10,
$ S := S^{(2,1)} S^{(1,2)}  \in L(L_p(0,1)$ factors as
$S = R^{(2)} R^{(1)} $ through $X $ with a dense range such that
each $R^{(i)}$ is compact.
Now
$$
T^{(1)} = R^{(1)} S^{(1,1)} \ , \ \ T^{(2)} = S^{(2,2)} R^{(2)}
$$
is the desired factorization. \quad\eop
\bigskip
{\bf 3.5. Corollary.} Every hilbertizable Fr\'echet--Schwartz space
has for all $1 < p < \infty $ and for all separable ${\cal L}_p$--spaces
$Y$ a system of local Banach spaces
isomorphic to  $Y$.
\bigskip
Let us still mention the following interesting case.
\bigskip
{\bf 3.6. Corollary.} {\sl Let  $(E,(p_{\alpha})_{\alpha \in A})$  be a
locally convex space such that every local Banach space  $E_{\alpha}
:= E_{p_{\alpha}}$  has the bounded approximation property and such that
for every  $\alpha$  there exist  $\beta,p_{\alpha} \leq p_{\beta}$,
and a complemented unconditional basic sequence
$(e_n(\alpha,\beta))_{n=1}^{\infty} \subset E_{\beta}$  such that the
restriction of the linking map  $T_{\alpha \beta} : E_{\beta} \to
E_{\alpha}$  to  ${\rm sp}\ (e_n(\alpha,\beta)\mid n \in {\Bbb N})$  is
compact.  Then  $E$  has a system of local Banach spaces
$(E_{q_{\alpha}})_{\alpha \in A}$  so that each  $E_{q_{\alpha}}$  has a
basis.}
\bigskip
{\bf Proof.}  It is easy to see, using Theorem 3.2, that each linking
map  $T_{\alpha \beta}$  as above factors densely through  $E_{\beta}
\times Z_{\alpha \beta}$  for an arbitrary separable Banach space
$Z_{\alpha \beta}$.  Since  $E_{\beta}$  has the bounded approximation
property, we can choose  $Z_{\alpha \beta}$  such that  $E_{\beta}
\times Z_{\alpha \beta}$  has a basis, see [P].  We then define the system of
seminorms as above Corollary 2.16.\quad\eop
\bigskip
{\bf Examples.} $1^{\circ}$  Every Schwartz space having a system of
local Banach spaces with the bounded approximation property and a
complemented unconditional basic sequence satisfies the conditions
of Corollary 3.6.

$2^{\circ}$.  Assume that  $E$  is a Fr\'echet space of Moscatelli
type (for definition, see [BD], Definition 1.3)
with respect to  $X, \ (Y_n)^{\infty}_{n=1}, \
(Z_n)^{\infty}_{n=1}$ and $(f_n)_{n=1}^{\infty}$,
where  $X$  is a normal Banach sequence space,
the Banach spaces  $Y_n$  and  $Z_n$  have the bounded approximation
property for all  $n$, the spaces  $Y_n$  have for all  $n$  unconditional
complemented basic sequences, and the linear maps $f_n: Y_n \to Z_n$ are
compact embeddings.  The conditions of
Corollary 3.6 are satisfied also in this case.
\bigskip
{\bf Remark.} The result that a Banach--space with an unconditional
basis is isomorphic to a complemented subspace of Banach space with a
symmetric basis (see [LT1], Theorem 3.b.1)
also yields an application similar to Corollary 3.6.

The following remark yields another method to combine
the results of Section 2
with the present considerations.
\bigskip
{\bf 3.7. Lemma.}  {\sl Let  $X,Y,Z$  and  $W$  be Banach spaces, $Z$
separable, let  $T \in L(X,Y)$  be compact and assume that every compact
operator  $S \in L(X,Y)$  factors densely through  $W$.  If  $X$  has
a complemented unconditional basic sequence  $(e_n)^{\infty}_{n=1}$,
then  $T$  factors densely through}  $W \times Z$.
\bigskip
{\bf Proof.} According to Theorem 3.2,  $T$  factors densely through
$X \times Z$  as  $T = S^{(2)}S^{(1)}$.  Moreover, the factorization
is such that  $S^{(2)}$  is compact.  By assumption, the restriction
of  $S^{(2)}$  to  $X$  factors densely through  $W,S^{(2)}|_X =
R^{(2)}R^{(1)}$, where  $R^{(1)} \in L(X,W),\ R^{(2)} \in L(W,X)$.  So,
writing
$$\eqalign{&T^{(1)} = (R^{(1)}P_X + P_Z)S^{(1)},\cr
           \noalign{\vskip 6pt}
           &T^{(2)} = R^{(2)}P_W + S^{(2)}P_Z,\cr}\eqno(3.5)$$
where  $P_X$  denotes the natural projection onto  $X$  etc.\ and
$T^{(1)} \in L(X,W \times Z),\ T^{(2)} \in L(W \times Z,X)$, we get
the desired factorization  $T = T^{(2)}T^{(1)}$.\quad\eop
\bigskip
{\bf Example.} Applying the results of Section 2 and Lemma 3.5 we see that
a compact operator on a separable reflexive Banach space  $X$  factors
densely through a separable Banach space  $W \times Z$, where  $W$  is
reflexive, if merely the uniform local distance of  $X^*$  and  $W^*$
is small enough and if some technical regularity
assumptions are satisfied.
To be more exact, the assumptions of Lemma 3.5 are satisfied, if  $X = Y$
and both  $X$  and  $W$  are reflexive, separable and have property
$(D)$  with constant functions  $f_X$  and  $f_W$, and if the uniform
local distance of  $X^*$  and  $W^*$  is bounded (Lemma 2.9).  For
example, every compact operator on a separable ${\cal L}_p$--space,
$1 < p < \infty$, factors densely through $\ell_p(E)$, the space of
$p$--summable, $E$  valued sequences, where  $E$  is an arbitrary
separable Banach space; see the example after Theorem 3.2.

Johnson constructed in [Jo] a family of separable Banach spaces  $C_p,\
1 \leq p \leq \infty$, with the strong property that every compact
operator on a separable Banach space  $X$  with the approximation
property factors through an  $C_p$.  The spaces $C_p$  are
$\ell_p$--sums of some sequences of finite dimensional Banach spaces.
Having a look at the proof of [Jo] one finds that the factorization
is not dense.  It would be interesting to know, if the factorization
could be done in a dense way for example in the case  $X$  has the
bounded approximation property.  Unfortunately, Theorems 2.8 or 2.14,
combined with Lemma 3.5, do not work here, since the local distance of
a space  $X$  and an $\ell_p$--sum of finite dimensional subspaces of
$X$  is usually not bounded.  The following result
implies a positive answer to this problem in a restricted case.
Note that because of the polar decomposition, the case of compact
operators in Hilbert space is  included there. For the properties
of symmetric bases we refer to [LT1], p. 113.
\bigskip
{\bf 3.8. Proposition.} {\sl Let  $X$  be a Banach space with a
normalized symmetric
basis  $(e_n)^{\infty}_{n=1}$  and let  $T \in L(X)$  be a compact
diagonal operator, $Te_n = \lambda e_n$  for a sequence
$(\lambda_n)^{\infty}_{n=1}$  of scalars satisfying  $\lim\limits_{n\to \infty}\
|\lambda_n| = 0$.  Let  $Y$  be a Banach space which has an
unconditional finite dimensional decomposition,
$ ( M_{n_k})_{k=1}^{\infty}$, where  $M_n = {\rm sp}\ (e_k| k \leq n)
\subset X$, for some increasing sequence  $(n_k)^{\infty}_{k=1}$.
If  $Z$  is an arbitrary separable Banach space,
then the operator  $T$  factors as  $T = T^{(2)}T^{(1)}$  through
$Y \times Z$  such that  $T^{(1)}(X)$  is dense in  $Y \times Z$  and
such that both  $T^{(1)}$  and  $T^{(2)}$  are compact.}
\bigskip
{\bf Proof.}  We may assume that
$|\lambda_n| \leq 1$  for all  $n$.  We choose the sequence
$(m_k)^{\infty}_{k=1},\ m_k \in {\Bbb N}$,  such that
$$|\lambda_{n}| \leq k^{-3}\eqno(3.6)$$
for all $n \geq m_k$,
for all  $k \in {\Bbb N}$.  Let
$$\eqalign{&X_1 = \overline{ {\rm sp}\ (e_n|n \neq m_k
\quad\hbox{for all}
\quad k) } \subset X ,\cr
\noalign{\vskip 6pt}
&X_2 = \overline{ {\rm sp}\ (e_{m_k}|k \in {\Bbb N})} \subset X.
\cr}\eqno(3.7)$$
Let  $({\nu}_k)^{\infty}_{k=1}$  be an increasing subsequence of
$(n_k)^{\infty}_{k=1}$  such that  ${\nu}_k \geq m_k$  for all  $k$.
We define
$$\eqalign{&Y_1 = \overline {\bigoplus_{k=1}^{\infty}\ M_{{\nu}_k}},\cr
\noalign{\vskip 6pt}
&Y_2 = \overline {\bigoplus_{k\in J}\ M_{n_k}}\cr}\eqno(3.8)$$
where  $J = \{k \in {\Bbb N}|n_k \neq {\nu}_t\quad\hbox{for all}\quad
t \in {\Bbb N}\}$.

The restriction of  $T$  to  $X_2$  factors densely through  $Y_2 \times
Z$  by (3.6) and Lemma 3.1.  Moreover, the factorization is such that
both factors are compact.  It is thus sufficient to prove that the
restriction of  $T$  to  $X_1$  factors densely through  $Y_1$  such
that both factors are compact.  We define the sets  ${\Bbb N}_k \subset
{\Bbb N},\ k = 1,2,\ldots$, inductively as follows: ${\Bbb N}_k$
consists of the smallest  ${\nu}_k$  natural numbers  $m$  such that
$m \neq m_t$  for all  $t$  and  $m \notin {\Bbb N}_t$  for  $t < k$.
Clearly, we have
$$m \geq {\nu}_{k-1} \geq m_{k-1}\eqno(3.9)$$
for all  $m \in {\Bbb N_k }$.  Moreover,
$$X_1 = \overline {\bigoplus_{k=1}^{\infty}\ N_k},\eqno(3.10)$$
where  $N_k = {\rm sp}\ \{e_n|n \in {\Bbb N}_k\}$.  Let  $\varphi_k:
\{1,\ldots,{\nu}_k\} \to {\Bbb N}_k$  be a bijection.  We define for all
$k \in {\Bbb N}$ the operator $ R^{(1)}_k \in L(N_k,M_{{\nu}_k})$  by
$$R^{(1)}_ke_{\varphi_k(n)} = |\lambda_{\varphi_k(n)}|^{1/2}e_n,\eqno(3.11)$$
where  $1 \leq n \leq {\nu}_k$,  and  $R_k^{(2)} \in L(M_{{\nu}_k},N_k)$  by
$$R_k^{(2)}e_n = \lambda_{\varphi_k(n)}|\lambda_{\varphi_k(n)}|^{-1/2}
e_{\varphi_k(n)}\eqno(3.12)$$
for  $1 \leq n \leq {\nu}_k$.  We have  $|\lambda_{\varphi_k(n)}| \leq
(k - 1)^{-3}$  for all  $k > 1 $  and for  $1 \leq n \leq {\nu}_k$, see the
choice of  $\varphi_k$, (3.9) and (3.6).  Since the basis
$(e_n)^{\infty}_{n=1}$  is symmetric, we get
$$\Vert R_k^{(1)}\Vert \leq Ck^{-3/2},\ \ \ \Vert R_k^{(2)}\Vert \leq
Ck^{-3/2}\eqno(3.13)$$
for a constant  $C$  (depending only on the properties of the basis
$(e_n)$; see [LT1], p.113)
and for all  $k$.  Denoting by  $(P_k)^{\infty}_{k=1}$  (resp.
$(Q_k)^{\infty}_{k=1}$) the uniformly bounded family of natural
projections, $P_k: X_1 \to N_k$  (resp.  $Q_k: Y_1 \to M_{{\nu}_k})$  we
define
$$R^{(1)} = \sum\limits_{k=1}^{\infty}\ R_k^{(1)}P_k,\ \
R^{(2)} = \sum\limits_{k=1}^{\infty}\ R_k^{(2)}Q_k.\eqno(3.14)$$
These operators are bounded and compact because of (3.13) and the
properties of  $(P_k)^{\infty}_{k=1}$  and  $(Q_k)^{\infty}_{k=1}$.
Moreover, (3.11) and (3.12) imply that  $R^{(2)}R^{(1)}e_n =
\lambda_ne_n$  for all  $n \in \cup_k\ {\Bbb N}_k$  so that
$R^{(2)} R^{(1)}x = Rx$  for  $x \in X_1$, see (3.10) and (3.7).  The
density of  $R^{(1)}(X_1)$  in  $Y_1$  also follows from
(3.11).\quad\eop
\bigskip
{\bf 3.9. Corollary.} {\sl Let  $X$  and  $T$  be as in Proposition 3.6.
For all  $p,\ 1 \leq p \leq \infty$, the operator  $T$  factors through
the Johnson space  $C_p$  such that the image of  $X$  in  $C_p$  is
dense.  As a consequence, each $\ell_q$--K\"othe sequence space,
$1 \leq q \leq \infty$,  which is a Schwartz space has a system of
local Banach spaces isomorphic to  $C_p$. }
\bigskip
{\bf Proof.} Given any sequence  $(M_n)^{\infty}_{n=1}$  of finite
dimensional Banach spaces the Johnson space contains a complemented
subspace  $Y$  isomorphic to  $(\oplus_{n=1}^{\infty}\ M_n)_{\ell_p}$.
So, Proposition 3.8 implies the desired factorization of  $T$.  The
statement concerning K\"othe spaces follows from the fact that the
linking maps between their natural local Banach spaces satisfy the
assumptions of Proposition 3.8.\quad\eop
\bigskip
We finally present a consequence or Theorem 3.2 containing Corollary 2.12.
\bigskip
{\bf 3.10 Corollary.} {\sl Let $Y$  be a Banach space. Every compact
operator $T \in L(\ell_2 , Y)$ factors densely through a separable
reflexive Banach space $X$ without basis.}
\bigskip
{\bf Proof.} Let $X$ be the Szarek space as in Proposition 2.4 with e.g.
$f(n) = 3n$. It is easy to see, using considerations like those
in the proof of Proposition 2.4, that for all $n \in {\Bbb N}$ we can
find $k$ such that $Y_k$ contains a $2$--complemented subspace
$2$--isomorphic to $\ell_{2}^{k}$. Since $X$ is the $\ell_2$--sum of the
spaces $Y_k$, we see that $X$ in fact contains a complemented copy of
$\ell_2$. Our result now follows directly from Theorem 3.2. \quad\eop
\bigskip
{\bf 4. Duality problems for local Banach spaces of Fr\'echet spaces.}
We consider the following problems on
the local Banach spaces of, say, a
Fr\'echet or a $(DF)$--space  $E$  and its strong dual  $E'_b$.
\bigskip
(L1) Assume that  $E$  has a system of local Banach spaces isomorphic to
a Banach space  $X$.  Does  $E'_b$  have a system of local Banach spaces
isomorphic to  $X^*$?

(L2) Is it possible to construct an example of a Fr\'echet or a
$(DF)$--space  $E$  such that, given large enough continuous
seminorms  $p$  and  $q$  on  $E$  and  $E'_b$, respectively, we have
$(E_\rho)^* \ncong (E'_b)_\gamma$  for all continuous seminorms  $\rho \geq p$
and  $\gamma \geq q$  on  $E$  and  $E'_b$?
\bigskip
There exists a simple counterexample to (L1). Let $E $ be any nuclear
$\ell_1$--K\"othe sequence
space with a continuous norm. Then $E $ is a Fr\'echet--space having
by definition  a system of local Banach spaces $(E_{p_k})_{k=1}^{\infty}$
isometric to $\ell_1.$ We have, for all $k \in {\Bbb N}$, $(E_{p_k})'
\cong \ell_{\infty}$, which is not a separable space. On the other hand
$E'_b$ is separable so that all the local Banach spaces of $E'_b$ are
also separable. So (L1) has a negative answer in this case.

The preceding counterexample is a separability argument, and it does
not give any information on the local structure of the Banach spaces
involved. In fact it is known that $E_b'$ has a system of local Banach
spaces isomorphic to $c_0$ (see [V]!). Both $c_0$ and $\ell_{\infty}$
are $ {\cal L}_{\infty}$--spaces so that at least in this sense they are still
quite similar. I do not think the above counterexample is yet a
satisfactory answer to problem (L1).

We can also ask the following natural question:
\bigskip
(B1) Assume that  $E$  has a system of local Banach spaces isomorphic to
the Banach space  $X$. Does $E$ also have
a family of Banach discs  $(B_{\alpha})_{\alpha \in A}$  such
that  $E_{B_{\alpha}} \cong X$  and such that every bounded set of  $E$
is contained in some  $B_{\alpha}$?
\bigskip
Clearly, questions (L1) and (B1) are related. To be more exact, every
weakly closed (with respect to the dual pair $<E,E'_b> $) absolutely
convex neighbourhood of zero $U \subset E'_b$ is the polar of a
Banach disc $B \subset E$, and, moreover, $(E')_U$ is isometric to
a subspace of $(E_B)'.$ The counterexample above does not
solve (B1).

We are not able to solve (B1) here.  I conjecture that the
answer is negative.  The difficulty in proving this
is that one has to be able to deal with a large class of Banach spaces
on which one has only little information and which may be very
pathological (cf.\ Corollary 2.17!).

To give some reference we mention the book [Ju] which contains a study of
related factorization problems.
The above problems are also connected with the problem of  projective
descriptions of  inductive limits. There exist a lot of papers on this topic by
K.D.Bierstedt, J.Bonet, R.Meise and others. We refer to the survey
articles [BM] and [BB].
In this context problem (B1) was solved in the positive
for $\ell_p$--K\"othe sequence spaces in [BMS], Proposition 2.5.

The result can easily be generalized to $X$--K\"othe sequence spaces in
the sense of Bellenot, [Be], where  $X$  is a Banach space with an
unconditional basis  $(e_n)^{\infty}_{n=1}$.  There is no difficulty to
give even a vector valued version.  Given  $X$  as above and a K\"othe
matrix  $(a_{kn})^{\infty}_{k,n=1}$  (i.e.\ a matrix consisting of non
negative numbers  $a_{kn}$  such that  $a_{k+1,n} \geq a_{kn}$  for all
$k,n$  and such that for all  $n$  there exists  $k$  with  $a_{kn} > 0)$
and a sequence  $(E_n,q_n)^{\infty}_{n=1}$  of Banach spaces, a vector
valued $X$--K\"othe sequence space  $E$  is defined by
$$E := \{x = (x_n)^{\infty}_{n=1}\mid x_n \in E_n,$$
$$p_k(x) := \Vert \sum_{n=1}^{\infty} a_{kn}q_n(x_n)e_n\Vert_X <
\infty\}.$$
Clearly, $E$  is a Fr\'echet space.  We assume now that  $a_{kn} > 0$
for all  $k$   and  $n$.  In this case all local Banach spaces  $E_{p_k}$
are isometric. Using the same method as [BMS], Proposition 2.5, we can now
prove
\bigskip
{\bf 4.1. Proposition.}  {\sl An arbitrary bounded set  $B \subset E$  is
contained in a bounded set of the form
$$B_0 = \{x\mid  \Vert \sum_{n=1}^{\infty} b_nq_n(x_n)e_n\Vert_X
\leq 1\}$$
for some positive sequence  $(b_n)^{\infty}_{n=1}$.  Moreover, $E_{B_0}$
is isometric to the spaces}  $E_{p_k}$.
\bigskip
{\bf Proof.} Let the numbers  $r_k,\ k \in {\Bbb N}$, be such that
$cB \subset r_kU_k$  and  $r_k < r_{k+1}$  for all  $k \in {\Bbb N}$,
where  $U_k := \{x \in E\mid p_k(x) \leq 1\}$ and
$c$ is the unconditionality constant of $(e_n).$  We define for all
$n \in {\Bbb N}$
$$b_n = \max\limits_{1\leq k\leq n} 2^{-k}r_k^{-1}a_{kn}.$$
Then the set  $B_0$, defined as above, is easily seen to be bounded.
We denote for all  $k \in {\Bbb N}$
$$K_k := \{n \in {\Bbb N} , \ n \geq k \mid
b_n = 2^{-k}r_k^{-1}a_{kn} , \ b_n \not= 2^{-t} r_t^{-1} a_{tn} \
{\rm for \ all } \ t,\  k < t \leq n \}.$$
Clearly, each $n \in {\Bbb N}$ belongs to exactly one $K_k$ so that
${\Bbb N}$ is a disjoint union of the sets $K_k$.
We have for all  $x = (x_n) \in \bigcap\limits_k\ r_kU_k$
$$\eqalign{&\Vert \sum_{n=1}^{\infty} b_nq_n(x_n)e_n\Vert_X\cr
\noalign{\vskip 8pt}
&= \Vert \sum_{k=1}^{\infty}\ \sum_{n\in K_k} 2^{-k}
r_k^{-1 } a_{kn}q_n(x_n)e_n\Vert_X\cr
\noalign{\vskip 8pt}
&\leq  \sum_{k=1}^{\infty} 2^{-k}r_k^{-1}\Vert\ \sum_{n\in K_k}
a_{kn}q_n(x_n)e_n\Vert_X\cr
\noalign{\vskip 8pt}
&\leq c \sup_{k\in {\Bbb N}}\ \{r_k^{-1}\Vert \sum_{n=1}^{\infty}
a_{kn}q_n(x_n)e_n\Vert_X. \} \cr}$$

This means that  $\bigcap\limits_{k=1}^{\infty} r_kU_k \subset cB_0$.

The last statement in our Proposition is clear from definitions.\quad\eop
\bigskip
So, in the case of vector valued $X$--K\"othe sequence spaces with
a continuous norm the answer to question (B1) is positive.  These spaces
are special examples of $T$--spaces defined in [BD] or (FG)--spaces
studied in [BDT].  It is an open problem if the answer to (B1) is
positive also in these more general classes.

Below we make an attempt to a negative solution of (B1).  This
construction shows the obstructions one has when trying to solve (B1)
in a positive direction.

We first concentrate on the following phenomenon.  Let  $A$  and  $B$
be closed absolutely convex sets in  ${\Bbb R}^n$  such that
${\rm sp}\ (A) = {\rm sp}\ (B) = {\rm sp}\ (A \cap B) = {\Bbb R}^n$,
and such that, say, $({\Bbb R}^n,A)$  and  $({\Bbb R}^n,B)$  are
isometric.  Then many isometric invariants occuring in Banach space
theory (like projection constants) may be very different for
$({\Bbb R}^n,A)$  and  $({\Bbb R}^n,A \cap B)$.  We give an example
which will be used to analyze problem (B1).

Note that if  $\Vert \cdot \Vert_A$  and  $\Vert \cdot \Vert_B$  are
the Minkowski functionals associated with  $A$  and  $B$, then the
Minkowski functional of $A \cap B$ equals
$x \mapsto \max\ \{ \Vert x\Vert_A,\Vert x\Vert_B\},\ x \in {\Bbb R}^n$.

We shall use tensor products and projection constants.  Let  $n \in
{\Bbb N},\ n = 2^k$  for some  $k \in {\Bbb N}$, and let  $(e_i)^n_{i=1}$
be the canonical basis of  $\ell^n_2$.  Let for all  $i,j,1\leq i,j \leq n$,
the numbers  $\varepsilon_{ij} = 1$  or  -1  be such that the matrix  $n^{-1/2}
(\varepsilon_{ij})^n_{i,j=1}$  is symmetric and orthogonal.  (See, for example,
[K1], 31.3.(5), p.429).  Let  $A: \ell^n_2 \to \ell^n_2$
be an operator such that
the matrix of  $A$  with respect to the basis  $(e_i)$  equals
$n^{-{1\over 2}}(\varepsilon_{ij})$.  Then also  $(Ae_i)_{i=1}^n$  is an
orthonormal basis of  $\ell_2^n$.

We denote by  $M$  the $n$--dimensional subspace of  $\ell^{2n}_2$,
spanned by the vectors  $f_i := (e_i,Ae_i),\ i = 1,\ldots,n$.  Note that
then also the vectors $g_i := (Ae_i,e_i)$  belong to  $M$, since
$$g_i = n^{-1/2}\sum_j\ \varepsilon_{ji}f_j\eqno(4.1)$$
because of the choice of the matrix  $n^{-1/2}(\varepsilon_{ij})$.

In the following we endow  $M$  with the norm of  $\ell_{\infty}^{2n}$,
which we denote by  $\Vert \cdot \Vert_{\infty}$.  The norm of
$(M,\Vert \cdot \Vert_{\infty})^*$  is denoted by
$\Vert \cdot \Vert_*$.

The following result is contained in [KTJ], Lemma 6. We give a different
proof.
\bigskip
{\bf 4.2. Lemma.} {\sl The absolute projection constant of $(M,
\Vert \cdot \Vert_{\infty})$  is at least}  $\sqrt n/2$.
\bigskip
{\bf Proof.}  Let  $(f^*_i)_{i=1}^n\ \ \hbox{(resp.}\ \  (g^*_i)^n_{i=1})
\subset M^*$  be the dual basis of  $(f_i)$ (resp. $ (g_i)$).
We consider the tensor
$$z = \sum_{i=1}^n f^*_i \otimes f_i \in M^* \otimes M \subset M^*
\otimes \ell_{\infty}^{2n}.$$
If  $z = \sum\ a_i \otimes b_i,\ a_i \in M^*,\ b_i \in M$  for all  $i$,
is an arbitrary finite representation, we have, by [Pi], B.1.4,
$$\sum_i\ \Vert a_i\Vert_*\Vert b_i\Vert_{\infty} \geq \sum_i\
\langle b_i,a_i\rangle = \sum_{i=1}^n \langle f_i,f_i^*\rangle = n.
\eqno(4.2)$$
Hence, the $M^* \otimes_{\pi} M$--norm of  $z$  is at least  $n$.

We claim that  $z$  has also the representation
$$z = \sum_{i=1}^n f_i^* \otimes (0,Ae_i) + \sum_{i=1}^n g_i^* \otimes
(Ae_i,0) \eqno(4.3)$$
in  $M^* \otimes \ell_{\infty}^{2n}$.  Indeed, the first sum in (4.3)
can be written as
$${1\over 2}\ \sum_i\ f_i^* \otimes ((e_i,Ae_i) + (-e_i,Ae_i)) \eqno(4.4)$$
and the second one as
$$\eqalignno{&{1\over 2}\ \sum_i\ g_i^* \otimes ((Ae_i,e_i) + (Ae_i,-e_i))\cr
\noalign{\vskip 8pt}
&= {1\over 2}\ \sum_i\ ({1\over \sqrt n}\ \sum_j\ \varepsilon_{ji}f^*_j) \otimes
({1\over \sqrt n}\ \sum_j\ \varepsilon_{ji}(e_j,Ae_j)\cr
\noalign{\vskip 8pt}
&\phantom{{1\over 2}\ \sum_i} + {1\over \sqrt n}\ \sum_j\
\varepsilon_{ji}
(e_j , - Ae_j))\cr
\noalign{\vskip 8pt}
&= {1\over 2} \ \sum_i\ f^*_i \otimes (e_i,A_i) + \sum_i\ f^*_i \otimes
(e_i,-Ae_i)&(4.5)\cr}$$
Note that here  $g^*_i = n^{-1/2}\sum\limits_j\ \varepsilon_{ji}f^*_j$  follows
from (4.1); moreover, for an orthonormal, symmetric $k\times k$--matrix
$(\sigma_{ij})^k_{i,j=1}$  we always have
$$\sum_{i=1}^k a_i \otimes b_i = \sum_{i=1}^k (\sum_{j=1}^k \sigma_{ji}a_j)
\otimes (\sum_{j=1}^k \sigma_{ji}b_j)$$
for vectors  $a_i,b_i$  in an arbitrary vector space.  So (4.3) follows
from (4.4) and (4.5).

We show that  $\Vert f_i^*\Vert,\Vert g_i^*\Vert \leq 1$  for all  $i$.
Let  $x \in M,\Vert x\Vert_{\infty} \leq 1$.  Then  $x$  has the
representations  \hbox{$x = \sum\limits_{i=1}^n a_if_i
= \sum\limits_{i=1}^n b_ig_i$},
and it follows from the definitions of  $f_i$  and  $g_i$  that
$\sup\ |a_i| \leq 1,\ \sup\ |b_i| \leq 1$.  Hence,
$$|\langle x,f^*_i\rangle | = |a_i| \leq 1,$$
and, similarly, $|\langle x,g^*_i\rangle | \leq 1$.  This proves that
$\Vert f^*_i\Vert,\Vert g^*_i\Vert \leq 1$.

Using (4.3) we now get the estimate
$$\eqalign{&\Vert z\Vert_{M'\otimes_{\pi}\ell_{\infty}^{2n}} \leq
\sum_{i=1}^n \Vert f^*_i\Vert \Vert (0,Ae_i)\Vert \cr
\noalign{\vskip 8pt}
&+ \sum_{i=1}^n \Vert g^*_i\Vert \Vert (Ae_i,0)\Vert \leq
2nn^{-1/2} = 2n^{1/2}.\cr}$$
This, combined with (4.2) implies  $\lambda (M) \geq n^{1/2}/2$  (see
for example [T], Lemma 4.1).\quad\eop
\bigskip
We now denote by  $\Vert \cdot \Vert_f$  and  $\Vert \cdot \Vert_g$
the $\ell^n_{\infty}$--norms in  $M$  with respect to the bases
$(f_i)^n_{i=1}$  and  $(g_i)^n_{i=1}$, respectively.  For
$1/\sqrt n < \alpha < \sqrt n$  we define the norm
$$\Vert x\Vert^{(\alpha)} := \max\ \{\Vert x\Vert_f,\alpha \Vert x\Vert_g\}
 \eqno(4.6)$$
in  $M$.
\bigskip
{\bf 4.3. Remark.}  The equality
$$\Vert x\Vert^{(1)} = \Vert x\Vert_{\infty}\eqno(4.7)$$
holds for all  $x \in M$.

Indeed, if  $x = \sum\limits_{i=1}^n x_if_i$, then we also have
$$\eqalignno{x &= \sum_{i=1}^n x_i\ \sum_{j=1}^n n^{-1/2}\varepsilon_{ji}g_j\cr
                  \noalign{\vskip 8pt}
               &= \sum_{j=1}^n (\sum_{i=1}^n n^{-1/2}x_i
               \varepsilon_{ji})g_j,&(4.8)\cr}$$
so that
$$\Vert x\Vert^{(1)} = \max\limits_{i,j}\ \{|x_i|,|\sum_{k=1}^n
n^{-1/2}x_k\varepsilon_{jk}|\}.\eqno(4.9)$$
On the other hand, if we consider  $x$  as an element of
$\ell_{\infty}^{2n}$,
$$\eqalign{&x = \sum_{i=1}^n x_i(e_i,Ae_i) = \sum_{i=1}^n x_i
(e_i,\sum_{j=1}^n n^{-{1\over 2}}\varepsilon_{ji}e_j)\cr
\noalign{\vskip 8pt}
&= (x_1,\ldots,x_n,\sum_{i=1}^n n^{-1/2}x_i\varepsilon_{1i},\ldots,\sum_{i=1}^n
n^{-{1\over 2}}x_i\varepsilon_{ni}).\cr}$$
This, together with (4.9), implies (4.7).
\bigskip
{\bf 4.4. Remarks.} We recall that every finite dimensional Banach space
$M$  is isometric to a subspace of  $c_0$.  Moreover, if  $X \subset c_0$
is a finite dimensional subspace, then the absolute projection constant
of $\lambda (X)$ satisfies
$$\lambda (X) = \lambda (X,c_0);\eqno(4.10)$$
see for example [TJ], Propositions 32.1 and 13.3.  We also remark that
if $X$  and  $Y$  are finite dimensional Banach spaces such that
$\lambda (X) > C$  and  $d(X,Y) \leq D$, then
$$\lambda (Y) > C/D.\eqno(4.11)$$
One can also easily prove that if  $Y$  is a Banach space isomorphic
to  $c_0,\
d(Y,c_0) < D \geq 1$  and  $X$  is a finite dimensional $C$--complemented
subspace of  $Y$, then
$$\lambda (X) \leq CD^2.\eqno(4.12)$$

The following lemma is essentially known; see for example [Ju],
Proposition 6.5.6 for the case $\alpha = 2 .$ For the sake of
completeness we give the proof.
\bigskip
{\bf 4.5. Lemma.} {\sl Let $(E, (p_k)^{\infty}_{k=1} )$ be a Fr\'echet space.
If $(r_k)_{k =1}^{\infty}$ is a sequence of positive numbers and $1 \leq
\alpha < \infty ,$ then the set
$$ B:= \{  x \in E \ \vert \ \left( \sum_{k=1}^{\infty} (r_k p_k (x)
)^{\alpha}  \right)^{1 / \alpha} \leq 1  \} $$
is a Banach disc. }
\bigskip
{\bf Proof.} It is easy to see that $B$ is bounded and
absolutely convex. That $B$ is closed can be seen as follows.
Let $\psi : E \to E_{p_k}$ be the canonical mapping induced by the
identity operator on $E$. Let us identify $E$ in the canonical way
with a subspace of $E_0 = \prod_{k =1}^{\infty} E_{p_k.}$ Then the
set $B $ is identified with
$$ B_0 = \{ x = (x_k )_{k=1}^{\infty}  \ \vert \ || x || :=
\left( \sum^{\infty }_{k=1} (r_k p_k (x_k) )^{\alpha} \right)^{1/ \alpha}
\leq 1, \ \exists y \in E \ {\rm such \ that \ } $$
$$ x_k = \psi_k y
\ {\rm for \ all \ } k \}.$$
Assume that $y = (y_k) \in \overline{B_0} \subset E \subset E_0 $
and that $(x^{(n)})_{n=1}^{\infty} \subset B_0$ is a sequence
converging to $y$ in the topology of $E_0.$ Since we have the product
topology on $E_0$, the sequence $(x^{(n)})$ converges coordinatewise
to $y.$ Given $\varepsilon > 0$ and $ m \in {\Bbb N}$ we thus find
$n \in {\Bbb N}$ such that
$$\sum^{m }_{k=1} (r_k p_k (x_k^{(n)}  - y_k ) )^{\alpha}
\leq \varepsilon^{\alpha}, $$
where  $ x^{(n)} = (x_k^{(n)})_{k}. $ Since $|| x^{(n)} || \leq 1, $
this implies
$ ( \sum^{m }_{k=1} (r_k p_k (y_k) )^{\alpha} )^{1/ \alpha}
<1 + \varepsilon ,$ and since $m$ is arbitrary, $|| y || \leq
1 + \varepsilon.$ Hence, $y \in B_0.$ \quad\eop
\bigskip
Using the remarks above we now construct an example of a
Fr\'echet space  $E$  which is a projective limit of Banach spaces
isomorphic to  $c_0$  such that  $E_B$  is not isomorphic to  $c_0$  for
''many'' bounded Banach discs  $B$.
\bigskip
{\bf 4.6. Proposition.}  {\sl There exists a separable Fr\'echet space
$(E,(U_k)^{\infty}_{k=1} ,(p_k)^{\infty}_{k=1})$  such that  $E_{p_k}
\cong c_0$  for all  $k$, and such that  $E_B \ncong c_0$, if  $B$ is
a Banach disc satisfying the following:}

(*) $p_B$  {\sl is formed using the real interpolation method from the norms
$p^{(\alpha)} :=$ \hfill \break$(\sum\limits_{k=1}^{\infty}
(r_k^{-1}p_k)^{\alpha})^{1/\alpha}$  and  $p^{(\beta)} :=
(\sum\limits_{k=1}^{\infty} (s_k^{-1}p_k)^{\beta})^{1/\beta}$  where
$1 \leq \alpha \leq \beta   \leq \infty$,
$(r_k)$  and  $(s_k)$  are arbitrary positive increasing sequences,
$r_k \leq s_k$  for all  $k$,
and}  $\sum\limits_{k=1}^{\infty} r_k^{-\alpha} \leq 1$.
\bigskip
To be more exact, we have $p_{\alpha} \geq p_{\beta}$ and thus the
Banach spaces $E_{\alpha}$ and $E_{\beta}$ corresponding to the
Banach discs $ \{ x \ \vert \ p^{(\alpha)} (x) \leq 1 \}$ and
$ \{ x \ \vert \ p^{(\beta)} (x) \leq 1 \},$ respectively, satisfy
$E_{\alpha } \subset E_{\beta}$ with a continuous embedding, and,
hence, they form a Banach interpolation couple (see [TJ], \S 3).
So, we can use the real interpolation method to produce norms $p_B$
defined on subspaces of $ E_{\alpha } + E_{\beta } = E_{\beta }.$
We always have $ p_{B } \geq p_{\beta}$ so that the closed unit ball
of $p_B$ is in fact a bounded disc. However, we do not claim that every
$B$ produced in such a way is closed in $E.$

Clearly, Proposition 4.6 contains as a special case the Banach discs
$$ B = \{ x \ \vert \ p(x) =
\left( \sum^{\infty }_{k=1} (r_k^{-1} p_k (x_k) )^{q} \right)^{1/ q} \leq 1
\}, $$
where $ 1 \leq q \leq \infty ;$ take $\alpha = \beta = q$ and $(s_k) =
(r_k)$ in Proposition 4.6. Note that if $q= \infty$, we have $B =
\bigcap_{k=1}^{\infty} r_k U_k.$

In fact, the only thing we need to assume on $B$ is that certain
projections on $E_B$ are well enough bounded. Unfortunately, it is hard
to describe this condition exactly before constructing the space $E.$
This is why we use the interpolation method to present our result.
We refer to the proof.
\bigskip
{\bf Proof.} We denote ${\Bbb N}_d = \{ 2^k | k \in {\Bbb N} \}$.
For all  $n  \in {\Bbb N}_d$, we denote by  $M_n$
the $n$--dimensional vector space  $M$  constructed above in this
section, and by  $(f_i^{(n)})^n_{i=1}$  and  $(g_i^{(n)})^n_{i=1}$  the
bases  $(f_i)^n_{i=1}$  and  $(g_i)^n_{i=1}$  of  $M$.  We also denote
by  $\Vert \cdot \Vert_{f,n}$  and  $\Vert \cdot \Vert_{g,n}$  the norms
$\Vert \cdot \Vert_f$  and  $n^{1/2}\Vert \cdot \Vert_g$  on  $M_n$.
Note that then  $\Vert x\Vert_{f,n} \leq \Vert x\Vert_{g,n}$  for all
$x \in M_n$, see (4.1) and the definition above (4.6).

Next we choose a bijection  ${\phi}: {\Bbb N}_d \to {\Bbb N} \times {\Bbb N}$
and define for all  $k \in {\Bbb N}$  the norms
$$ p_k((x_n)_{n \in  {\Bbb N}_d } ) = \sup_{n\in {\Bbb N}_d} \
\{ p_{k,n} (x_n) \}
\eqno(4.13)$$
where  $(x_n) \in \bigoplus_{{\Bbb N}_d} M_n$  and
$$p_{k,n}(x_n) := \cases{\Vert x_n\Vert_{f,n}, &if
$\pi_1(\phi (n)) \geq k$\cr
\noalign{\vskip 6pt}
\Vert x_n\Vert_{g,n}, &if  $\pi_1(\phi ( n)) < k$;\cr}\eqno(4.14)$$
here  $\pi_i,\ i = 1,2$, denote the canonical projections of  ${\Bbb N}
\times {\Bbb N}$  onto the first and second coordinate spaces,
respectively.  Taking the completion of  $\bigoplus_{n \in {\Bbb N}_d}
M_n$  with
respect to the topology determined by the norms  $p_k$  we get a
separable Fr\'echet space  $(E,(U_k)^{\infty}_{k=1},(p_k)^{\infty}_{k=1})$.
Clearly, each local Banach space  $E_{p_k}$  is isometric to  $c_0$.

Assume now  that  $B \subset E$  is a Banach disc satisfying $(*)$ for
some $\alpha $ and $\beta $ and for some sequences  $(r_k)$  and  $(s_k)$.  We show that  $E_B$  has an
infinite family  $(X_m)^{\infty}_{m=1}$  of finite dimensional
1--complemented subspaces such that
$$\sup_{m\in {\Bbb N}} \{ \lambda (X_m) \} = \infty.\eqno(4.15)$$
We choose for all  $m \in {\Bbb N}$  the number  $K$  such that
$(\sum\limits_{k=K}^{\infty} r_k^{-\alpha})^{-1/\alpha} > 4 ms_1$,
and then  $N \in {\Bbb N}_d$
such that  $\pi_1(\phi  (N)) = K - 1$  and $N^{1/2} >
( \sum\limits_{k=K}^{\infty} r_k^{-\alpha})^{-1/\alpha} $ and
$N > s_1^{-2} s^2_K$.  Let  $P$
be the natural projection from  $\bigoplus\limits_n\ M_n$  onto  $M_N$.
Now  $P$  is continuous with the operator norm equal to  1  when
$\bigoplus\limits_n\ M_n$  is endowed with either of the norms
$p^{(\alpha)}$  or  $p^{(\beta)}$.  Since  $p_B$  is formed by
interpolation, we see that  $\Vert P\Vert = 1$  also as an operator
in  $E_B$  (cf.\ [TJ], \S 3), so that  $M_N$  is a 1--complemented
subspace of  $E_B$. (This is the only point where we need to use  some
extra assumption on $B.$)

Using the properties of the sequence  $(r_k)$  and the definition of
the norms  $p_k$  we get
$$\eqalignno{&p^{(\alpha)}\big|_{M_N} = (\sum^{K-1}_{k=1} (r_k^{-1}p_k
\big|_{M_N})^{\alpha} + \sum_{k=K}^{\infty}
(r_k^{-1}p_k\big|_{M_N})^{\alpha})^{1/\alpha}\cr
\noalign{\vskip 8pt}
&\leq (\sum^{K-1}_{k=1} (r_k^{-1} || \cdot ||_{f,N} )^{\alpha}
+ \sum_{k=K}^{\infty}
(r_k^{-1} || \cdot ||_{g,N} )^{\alpha})^{1/\alpha}
\cr
\noalign{\vskip 8pt}
&\leq (R_1^{\alpha}( || \cdot ||_{f,N} )^{\alpha} +
R_K^{\alpha}(|| \cdot ||_{g,N} )^{\alpha})^{1/\alpha}\cr
\noalign{\vskip 8pt}
&\leq 2\max\ \{R_1\Vert \cdot \Vert_{f,N},R_K\Vert \cdot \Vert_{g,N}\},
&(4.16)\cr}$$
where  $R_{k'} = (\sum\limits_{k=k'}^{\infty} r_k^{-\alpha})^{1/\alpha}$.
On the other hand, obviously
$$\max\ \{ s_1^{-1}\Vert \cdot \Vert_{f,N},s_K^{-1}\Vert
\cdot \Vert_{g,N} \}
\leq p^{(\beta) }\big|_{M_N}.  \eqno(4.17)$$
By the properties of interpolation, we have $ p_{\alpha }\big|_{M_N}
\geq p_B \big|_{M_N}    \geq p_{\beta} \big|_{M_N} .$
Combining this with (4.16) and (4.17) and the assumption
$R_1 \leq 1$ (see (*)) yields
$$\eqalign{&2 \max\ \{ \Vert x\Vert_{f,N},R_K\Vert x\Vert_{g,N}\}\cr
\noalign{\vskip 6pt}
&\geq p_B(x) \geq \max\ \{s_1^{-1}\Vert x\Vert_{f,N},s_K^{-1}
\Vert x\Vert_{g,N}\}\cr}$$
for  $x \in M_N$, and, using the notation (4.6) and the definitions of
$\Vert \cdot \Vert_{f,N}$  and  $\Vert \cdot \Vert_{g,N}$, we get
$$\eqalignno{2 R_KN^{1/2}\Vert x\Vert^{(1)} & \geq
2 \max \{ ||x ||_{f,N} , \ R_K ||x ||_{g,N} \} \geq  p_B(x)\cr
\noalign{\vskip 8pt}
&\geq  s_1^{-1}\Vert x\Vert^{(N^{1/2}s_1s_K^{-1})} \geq
s_1^{-1}\Vert x\Vert^{(1)}& \cr}$$
for all  $x \in M_N$, by the choice of  $N$.  Hence,
$$d((M_N,p_B),(M_N,\Vert \cdot \Vert^{(1)})) \leq 2  s_1 R_KN^{1/2}.
\eqno(4.18)$$
On the other hand, by Lemma 4.2 and Remark 4.3,  $\lambda ((M_N,
\Vert \cdot \Vert^{(1)})) \geq \sqrt N/2$.  So, by the remark (4.11) and
the choice of  $K$,
$$\lambda ((M_{N},p_B)) \geq  s_1^{-1}R_K^{-1}/ 4 \geq m,$$
which proves (4.15).

By the Remark 4.4 it is now clear that  $E_B$  cannot be isomorphic to
$c_0$, since  $E_B$  contains 1--complemented subspaces, the absolute
projection constants of which can be chosen arbitrarily large.\quad\eop
\bigskip
{\bf 4.7. Remark.} The last conclusion in the proof of Proposition 4.6 even
shows that the local distance of $c_0$ and $E_B$ is not bounded. Hence,
the counterexample given in Proposition 4.6. is much stronger than the
trivial counterexample in the beginning of this section.
\bigskip\bigskip\bigskip
\centerline{\bf References}
\bigskip
\ref [BB] Bierstedt, K.D., Bonet, J.: Projective descriptions of
weighted inductive limits: the vector valued cases. Editor T.Terzioglu,
NATO ASI Series C, Vol 287, Kluwer Academic Publishers (1989), 195--221.

\ref [BM] Bierstedt, K.D., Meise, R.: Weighted inductive limits and
their projective descriptions. DOGA Tr.J.Math. 10,1 (1986), 54--82.

\ref [BMS] Bierstedt, K.D., Meise, R., Summers, W.: K\"othe sets
and K\"othe sequence spaces.  Functional Analysis, Holomorphy and
Approximation Theory.  J.A.\ Barroso (ed.) North Holland Math.\ Studies
71 (1982) p.\ 27--91.

\ref [BDG] Bonet, J., Defant, A., Galbis, A.: A note on two questions
of Grothendieck and their relation to the geometry of Banach spaces.
Preprint.

\ref [BD] Bonet, J., Diaz, J.C.: The problem of topologies of
Grothendieck and the class of Fr\'echet T--spaces. Math. Nachr. 150
(1991), 109--118.

\ref [BDT] Bonet, J., Diaz, J.C., Taskinen, J.: Tensor stable Fr\'echet
and (DF)--spaces. Preprint.

\ref [BDi] Bonet, J., Dierolf, S.: Fr\'echet spaces of Moscatelli type.
Rev. Mat. Univ. Complut. Madrid 2 (suppl.)(1989), 77--92.

\ref [DFJP] Davis, W.J., Figiel, T., Johnson, W.B., Pe{\l}czy\'nski, A.:
Factoring weakly compact operators. J.Funct. Anal. 17 (1974), 311--327.

\ref [H] Hogbe--Nlend, H.: Les espaces de Fr\'echet--Schwartz et la
propri\'et\'e d'approximation. C.R.Acad.Sci.Paris (A) 275 (1972),
1073--1075.

\ref [J] Jarchow, H.: Locally convex spaces.  Teubner, Stuttgart (1981).

\ref [Jo] Johnson, W.B.: Factoring compact operators. Israel J. Math. 9
(1971), 337--345.

\ref [JP] Johnson, W.B., Pisier, G.: The proportional UAP characterizes
weak Hilbert spaces. Preprint.

\ref [JRZ] Johnson, W.B., Rosenthal, H.P.,  Zippin, M.: On bases, finite
dimensional decompositions and weaker structures in Banach spaces.
Israel J.\ Math.\ 9 (1971), 488--506.

\ref [Ju] Junek, H.: Locally convex spaces and operator ideals. Teubner,
Leipzig, (1983).

\ref [KP] Kadec, M.I, Pe{\l}czy\'nski, A.: Bases, lacunary sequences and
complemented subspaces in the spaces $L_p$. Studia Math. 21 (1962),
161--176.

\ref [KTJ] K\"onig, H., Tomczak--Jaegermann, N.: Bounds for projection
constants and 1--summing norms. Preprint.

\ref [K1] K\"othe, G.: Topological vector spaces I, 2.\ ed.\ Springer,
Berlin--Heidelberg--New York (1983).

\ref [K2] K\"othe, G.: Topological vector spaces II,  Springer,
Berlin--Heidelberg--New York (1979).

\ref [L] Lewis, D.: Finite dimensional subspaces of  $L_p$.  Studia
Math., 63 (1978), 207--212.

\ref [LP] Lindenstrauss, J., Pe{\l}czy\'nski, A.: Absolutely summing
operators  in ${\cal L}_p$--spaces and their applications. Studia Math.
29 (1968), 275--326.

\ref [LT] Lindenstrauss, J., Tzafriri, L.: Classical Banach Spaces.
Springer Lecture Notes in Mathematics 338 (1973).

\ref [LT1] Lindenstrauss, J., Tzafriri, L.: Classical Banach Spaces
I.  Springer, Berlin--Heidel\-berg--New York (1977)

\ref [M] Mascioni, V.: On the duality of the uniform approximation
property in Banach spaces. Preprint.

\ref [MS] Milman, V., Schechtman, G.: Asymptotic
theory of finite dimensional normed spaces.
Springer Lecture Notes in Mathematics 1200 (1986).

\ref [P] Pe{\l}czy\'nski, A.: Any separable Banach space with the
bounded approximation property is a complemented subspace of a Banach
space with a basis. Studia Math. 40 (1971), 239--243.

\ref [PR] Pe{\l}czy\'nski, A., Rosenthal, H.P.: Localization techniques
in $L^p$--spaces. Studia Math. 52 (1975), 263--289.

\ref [Pi] Pietsch, A.: Operator Ideals.  North Holland (1978).

\ref [Ps1] Pisier, G.: Counterexamples to a conjecture of Grothendieck.
Acta Math.\ 151 (1983), 181--208.

\ref [Ps2] Pisier, G.: Weak Hilbert spaces. Proc. London Math. Soc.
(3) 56 (1988), 547--579.

\ref [S] Szarek, S.: A Banach space without a basis which has the
bounded approximation property.  Acta Math.\ 159 (1987), 81--98.

\ref [T] Taskinen, J.: Counterexamples to ''Probl\`eme des topologies''
of Grothendieck.  Ann.\ Acad.\ Sci.\ Fenn.\ Ser.\ A I Math.\ Diss.\ 63
(1986).

\ref [TJ] Tomczak--Jaegermann, N.: Banach--Mazur distances and
finite--dimensional operator ideals.  Pitman Monographs and Surveys
in Pure and Applied Mathematics 38.

\ref [V] Valdivia, M.: Nuclearity and Banach spaces. Proc.Edinburgh
Math. Soc. 20 (1977), 205--209.
\bigskip\bigskip\bigskip
\noindent Department of Mathematics\hfill\break
\noindent University of Helsinki\hfill\break
\noindent Hallituskatu 15\hfill\break
\noindent 00100 HELSINKI\hfill\break
\noindent Finland

\bye